\documentclass[11pt]{amsart}
\usepackage{a4wide,amssymb}
\usepackage{todonotes}
\usepackage{hyperref}
\date{}
\newtheorem{theorem}{Theorem}[section]
\newtheorem{lemma}[theorem]{Lemma}
\newtheorem{corollary}[theorem]{Corollary}
\newtheorem{proposition}[theorem]{Proposition}
\newtheorem{claim}[theorem]{Claim}
\newtheorem{question}[theorem]{Question}
\newtheorem{problem}[theorem]{Problem}
\theoremstyle{definition}

\newtheorem{remark}[theorem]{Remark}
\newtheorem{example}[theorem]{Example}

\def\cc {{\mathfrak c}}
\newcommand{\0}{\emptyset}
\newcommand{\U}{\mathcal U}
\newcommand{\V}{\mathcal V}
\def\phi{\varphi}

\def\bs{\backslash}
\def\ol{\overline}
\def\RT{\Rightarrow}

\def\supp{\operatorname{supp}}

\def\int{\operatorname{int}}
\newcommand{\w}{\omega}
\newcommand{\Ra}{\Rightarrow}
\newcommand{\IR}{\mathbb R}
\newcommand{\IZ}{\mathbb Z}
\newcommand{\IQ}{\mathbb Q}
\newcommand{\IC}{\mathbb C}
\newcommand{\IT}{\mathbb T}

\begin{document}
\title[On feebly compact paratopological groups]
{On feebly compact paratopological groups}
\author[Taras Banakh]{Taras Banakh}
\address{Taras Banakh: Ivan Franko National University of Lviv (Ukraine), and Jan Kochanowski University in Kielce (Poland)}
\email{t.o.banakh@gmail.com}
\author[Alex Ravsky]{Alex Ravsky}
\address{Alex Ravsky: Department of Analysis, Geometry and Topology, Pidstryhach Institute for Applied Problems of Mechanics and Mathematics
National Academy of Sciences of Ukraine,
Naukova 3-b, Lviv, 79060, Ukraine}
\email{alexander.ravsky@uni-wuerzburg.de}
\begin{abstract}
We obtain many results and solve some problems about feebly compact paratopological groups.
We obtain necessary and sufficient conditions for such a group to be topological. One of them is the
quasiregularity. We prove that each $2$-pseudocompact paratopological group is feebly compact and that
each Hausdorff $\sigma$-compact feebly compact paratopological group is a compact topological group.
Our particular attention concerns periodic and topologically periodic groups.
We construct examples of various compact-like paratopological groups which are not topological
groups, among them a $T_0$ sequentially compact group, a $T_1$ $2$-pseudocompact group, a
functionally Hausdorff countably compact group (under the axiomatic assumption that there is an
infinite torsion-free Abelian countably compact topological group without non-trivial convergent
sequences), and a functionally Hausdorff second countable group sequentially pracompact group.
We prove that the product of a family of feebly compact paratopological
groups is feebly compact, and that a paratopological group $G$ is feebly compact provided
it has a feebly compact normal subgroup $H$ such that a quotient group $G/H$ is feebly compact. For our research we also study some general constructions of paratopological groups.
We extend the well-known construction of Ra\u\i kov completion of a $T_0$ topological group
to the class of paratopological groups.
We investigate cone topologies of paratopological groups
which provide a general tool for constructing pathological examples,
especially examples of compact-like paratopological groups with discontinuous inversion.
We find a simple interplay between the algebraic properties
of a basic cone subsemigroup $S$ of a group $G$ and compact-like properties of two basic semigroup
topologies generated by $S$ on the group $G$.
\end{abstract}

\maketitle 
\section*{Introduction}
In this paper we study paratopological groups that are feebly compact or have other compact-like properties. The investigation was motived by the classical problem of finding compact-like properties under which a paratopological group is a topological group. This problem is discussed in Section~\ref{PPGTAT}. In Section~\ref{s:Raikoff} we extend the well-known construction of
Ra\u\i kov completion of a $T_0$ topological group to the class of paratopological groups.
In Section~\ref{CCSProd} we use the results of Sections~ \ref{s:Raikoff} and \ref{PPGTAT} to prove that the feeble compactness is preserved by Tychonoff products and extensions of paratopological groups. In Section~\ref{sec:cone-top} we analyse the so-called cone topologies on Abelian groups, which allows us to construct many ``pathological'' examples of paratopological groups distinguishing between various compact-like properties. More detail information on main results and (un)solved problems can be found in the introductions of the sections.
\section{Preliminaries}
In this paper the word ``space" means ``topological space''.
\subsection{Topologized groups}
A pair $(G,\tau)$ consisting of a group $G$ and a topology $\tau$ on $G$ is called a {\em topologized group}.
A topologized group $(G,\tau)$ is called
\begin{itemize}
\item a {\em semitopological group} if the topology $\tau$ is {\em shift-continuous}, that is for any $a,b\in G$ the two-sided shift $s_{a,b}:G\to G$, $s_{a,b}:x\mapsto axb$, is continuous;
\item a {\em quasitopological group} if the topology $\tau$ is shift-continuous and the inversion $G\to G$, $x\mapsto x^{-1}$, is continuous;
\item a {\em paratopological group} if $\tau$ is a {\em semigroup topology}, which means that the multiplication $G\times G\to G$, $(x,y)\mapsto xy$, is continuous;
\item a {\em topological group} if $\tau$ is a {\em group topology} which means that $\tau$ is a semigroup topology with continuous inversion $G\to G$, $x\mapsto x^{-1}$.
\end{itemize}
It follows that a topologized group is a topological group if and only if it is both a paratopological group and a quasitopological group.
For a shift-invariant topology $\tau$ on a group $G$ by $\tau^\sharp$ we denote the shift-invariant topology on $G$ whose neighborhood base at the identity $e$ of $G$ consists of the sets $U\cap U^{-1}$ where $e\in U\in\tau$. It follows that $(G,\tau^\sharp)$ is a quasitopological group, which is topologically isomorphic to the subgroup $\{(x,x^{-1}):x\in G\}$ of  $G\times G$. If $(G,\tau)$ is a paratopological group, then $(G,\tau^\sharp)$ is a topological group, called the {\em group coreflexion} of $(G,\tau)$.
In this paper we shall be mainly interested in paratopological groups.
A classical example of a paratopological group failing to be a
topological group is the Sorgenfrey line, that is the real
line endowed with the Sorgenfrey topology (generated by the
base consisting of half-intervals $[a,b)$, $a<b$).
For a semitopological group $(G,\tau)$ by $\tau_e=\{U\in\tau:e\in U\}$ we denote the family of $\tau$-open neighborhoods of the identity $e$ of the group.

Whereas the investigation of topological groups already is a fundamental branch of topological
algebra (see, for instance,~\cite{Pon},~\cite{DikProSto}, and~\cite{ArhTka}),
semitopological groups are not so well-studied and have more variable structure.
Basic properties of semitopological or paratopological groups compared with
the properties of topological groups are described in the book \cite{ArhTka} by Arhangel'skii and Tkachenko, in the PhD thesis of Ravsky~\cite{Rav3},
papers~\cite{Rav}, ~\cite{Rav2}, and in the survey \cite{Tka5} by Tkachenko.
\subsection{Separation axioms} These axioms describe specific structural properties of spaces.
Basic separation axioms and relations between them are considered in~\cite[Section 1.5]{Eng}.
For more specific cases and topics, also related to semitopological and paratopological groups,
see ~\cite[Section~1]{Rav2}, ~\cite{BanRav2},
~\cite[Section~2]{Tka5},~\cite{Tka7},~\cite{LiTuXie}.
All spaces considered in the present paper are {\it not} supposed to satisfy any of the separation
axioms, if otherwise is not stated. For a subset $A$ of a space $(X,\tau)$ by $\overline{A}^\tau$
and $\mathrm{int}_\tau A$ we denote the closure and interior of $A$ in $X$, respectively. If the
topology $\tau$ is clear from the context, then we write $\overline{A}$ and $\mathrm{int} A $
instead of $\overline{A}^\tau$ and $\mathrm{int}_\tau(A)$, respectively. Also by
$\tau{\restriction}A$ we denote the subspace topology $\{U\cap A:U\in\tau\}$ on $A$, inherited from
the topology $\tau$ of $X$.
\smallskip
Now we recall the definition of some separation axioms. A space $X$ is
\begin{itemize}
\item {\it $T_0$} if for any distinct points $x,y\in X$ there exists an open set
$U\subseteq X$, which contains exactly one of the points $x$, $y$;
\item {\it $T_1$} if for any distinct points $x,y\in X$ there exists a open set
$U\subseteq X$ such that $x\in U$ and $y\notin U$;
\item {\it $T_2$} or {\it Hausdorff\/} if any distinct points $x,y\in X$ have disjoint
neighborhoods in $X$;
\item {\it regular} if any neighborhood of any point contains a closed neighborhood of the point;
\item {\it $T_3$} if it is regular $T_0$;
\item {\it quasiregular} if any nonempty open subset of $X$ contains the closure
of some nonempty open subset of $X$;
\item {\it semiregular} if $X$ has a base consisting of {\it regular open} sets,
that is sets $U=\int\ol U$;
\item {\it functionally Hausdorff\/} if for any distinct points $x,y\in X$
there exists a continuous function $f:X\to\mathbb R$ such that $f(x)\ne f(y)$;
\item {\it completely regular} if for any closed set
$F\subseteq X$ and  point $x\in X\setminus F$ there exists a continuous function
$f:X\to\mathbb R$ such that $f(x)=0$ and $f(F)\subseteq \{1\}$;
\item $T_{3\frac12}$ or {\em Tychonoff\/} if it is completely regular  and $T_0$.
\end{itemize}
Remark that each regular space is quasiregular and semiregular. 
If a space $(X,\tau)$ satisfies a separation axiom $T$, we shall say that a
topology $\tau$ is $T$. For instance, if $(X,\tau)$ is regular then $\tau$ is a regular topology.
Given a space $(X,\tau)$, Stone \cite{Sto} and Kat\v{e}tov \cite{Kat} considered the topology
$\tau_{sr}$ on $X$ generated by the base consisting of all regular open sets of the space
$(X,\tau)$. The topology $\tau_{sr}$ is called the {\it semiregularization} of the topology $\tau$
and the space $(X,\tau_{sr})$ is called the {\it semiregularization} of the space $(X,\tau)$.
It is easy to see that the semiregularization of a Hausdorff space is Hausdorff.
\subsection{Separation axioms in paratopological groups}
It is well-known that  each topological group is completely regular.
On the other hand, simple examples show that none of the implications
$T_0\Rightarrow T_1 \Rightarrow T_2 \Rightarrow T_3$ holds for paratopological groups (see \cite[Examples 1.6-1.8]{Rav} and
page 5 in any of papers \cite{Rav2} or \cite{Tka5})
and there are only a few backwards implications between different separation axioms, see
~\cite[Section 1]{Rav2} or~\cite[Section 2]{Tka5}.
On the other hand, the authors proved in \cite{BanRav2} that every semiregular paratopological group is completely regular and every Hausdorff paratopological group is functionally Hausdorff. The proof of the following important fact can be found in  \cite[Ex. 1.9]{Rav2}, \cite[p. 31]{Rav3} or \cite[p. 28]{Rav3}.
\begin{lemma}\label{l:semiregularization} For any (Hausdorff) paratopological group $(G,\tau)$,
its semiregularization $(G,\tau_{sr})$ is a (Hausdorff) regular paratopological group.
\end{lemma}
Finally, let us mention that the paper~\cite{Tka7} discusses some functors related to  separation axioms in paratopological groups.


\subsection{Compact-like spaces} Different classes of compact-like spaces and relations between them
provide a well-known investigation topic in General Topology, see, for instance, basic \cite[Chap.
3]{Eng} and general \cite{Vau}, \cite{Ste}, \cite{DRRT}, \cite{Mat}, \cite{Lip} works. The
inclusion relations between the classes are often visually represented by arrow diagrams, see,
\cite[Diag. 3 at p.17]{Mat}, \cite[Diag. 1 at p. 58]{Dor-AldSha} (for completely regular spaces),
\cite[Diag. 3.6 at p. 611]{Ste}, and~\cite[Diag. at p. 3]{GutRav}. For recent results on
embeddings of spaces into various Hausdorff compact-like spaces, see \cite{BanM}, \cite{BBR1},
\cite{BBR2}.
Let us recall the definitions of compact-like spaces which will appear in this paper.
A space $X$ is called
\begin{itemize}
\item {\it sequentially compact} if each sequence in $X$ contains
a convergent subsequence;
\item {\it $\omega$-bounded} if each countable subset of $X$ has compact closure;
\item {\it totally countably compact} if each sequence of $X$ contains
a subsequence with  compact closure;
\item {\it countably compact at a subset} $A$ of $X$, if each infinite
subset $B$ of $A$ has an accumulation point $x\in X$
(the latter means that each neighborhood of $x$ contains infinitely many points of the set $B$);
\item {\it countably compact} if $X$ is countably compact at itself, or, equivalently,
if each countable open cover of $X$ has a finite subcover;
\item {\it sequentially pracompact} if $X$ contains a dense set $D$ such that
each sequence of points of the set $D$ has a convergent subsequence;
\item {\it countably pracompact} if $X$ is countably compact at a dense subset of $X$;
\item {\it feebly compact} if each locally finite family of open sets of the space $X$ is
finite;
\item {\it sequentially feebly compact} if for every sequence $(U_n)_{n\in\w}$ of nonempty open sets of $X$ there exists a point $x\in X$ and an infinite set $I\subseteq\w$ such that for every neighborhood $O_x\subseteq X$ of $x$ the set $\{n\in I:O_x\cap U_n\ne\emptyset\}$ is finite;
\item {\it pseudocompact} if $X$ is Tychonoff and each continuous real-valued function on $X$ is bounded;
\item {\it Lindel\"of} if each open cover of $X$ has a countable subcover;
\item {\em $H$-compact} if any open cover $\U$ of $X$ contains a finite subfamily $\V\subseteq\U$ such that $X=\bigcup_{V\in\V}\overline{V}$;
\item {\em $H$-closed} if $X$ is $H$-compact and Hausdorff;
\item {\it countably o-compact} if each decreasing sequence $(U_n)_{n\in\omega}$ of nonempty
open sets in $X$ has nonempty intersection $\bigcap_{n\in\omega}U_n$.
\end{itemize}
A sequence $\{U_n:n\in\omega\}$ of subsets of a space $X$ is {\it decreasing} if
$U_n\supseteq U_{n+1}$ for every $n\in\omega$.

These notions relate as follows:
\begin{itemize}
\item Each compact space is $\omega$-bounded and $H$-compact.
\item Each $\omega$-bounded space is totally countably compact.
\item Each totally countably compact space is countably compact.
\item Each sequentially compact space is countably compact and  sequentially pracompact.
\item Each countably compact space is countably pracompact.
\item Each sequentially pracompact space is countably pracompact.
\item Each countably pracompact space is feebly compact.
\item Each countably o-compact space is feebly compact.
\item Each Hausdorff countably o-compact space is finite.
\item A space is compact if and only if it is countably compact and Lindel\"of.
\item A space is pseudocompact if and only if it is Tychonoff and feebly compact, see \cite[Theorem 3.10.22]{Eng}.
\item  A space is $H$-compact if and only if its semiregularization is $H$-compact.
\item A Hausdorff space $X$ is $H$-compact if and only if $X$ is closed in any Hausdorff space containing $X$ as a closed subspace, see \cite[3.12.5]{Eng}.
 \end{itemize}
\subsection{Compact-like semitopological groups}
A semitopological group $G$ is {\it left} (resp. {\em right}) {\em precompact } if for each neighborhood $U$ of
the identity of $G$, there exists a finite subset $F$ of $G$ such that $FU=G$ (resp. $UF=G$).
Proposition 2.1 from~\cite{Rav2} implies that a
paratopological group is left precompact iff it is right precompact.
So we will call left precompact paratopological groups {\it precompact}.
 A semitopological group $G$ is {\it $2$-pseudocompact}
if $\bigcap_{n\in\omega} \ol{U_n^{-1}}\not=\0$ for any decreasing sequence $\{U_n:n\in\omega\}$ of nonempty
open subsets of $G$. It is easy to see that a semitopological group is $2$-pseudocompact if it
is countably compact or countably o-compact.  In Proposition~\ref{2PseP} we shall prove that each $2$-pseudocompact paratopological group is feebly compact.

A semitopological group $G$ is {\it left
$\omega$-precompact} if for each neighborhood $U$ of the identity of $G$ there exists a countable
subset $F$ of $G$ such that $FU=G$.
A semitopological group $G$ is {\em saturated} if for each nonempty open subset $U$ of $G$
there exists a nonempty open subset $V$ of $G$ such that $V^{-1}\subseteq U$.
By~\cite[Proposition 3.1]{Rav2} each precompact paratopological group is saturated
and by~\cite[Proposition 1.7]{Rav2} each saturated paratopological group is quasiregular.
Recent investigations in papers ~\cite{Aza}, ~\cite{LinLin},
and \cite{LinLinSan} are using the following notions. A topologized group $G$ is called
{\it pseudobounded} provided for each neighborhood $U$ of the identity of the group $G$
there exists a natural number $n$ such that $U^n=G$.
A topologized group $G$ is called
{\it $\omega$-pseudobounded} provided $G=\bigcup_{n\in\mathbb N} U^n$ for any neighborhood $U$ of the identity in $G$.
Clearly, every pseudobounded topologized group is $\omega$-pseudobounded. However, the additive
group $(\mathbb R, +)$ endowed with the standard topology is an $\omega$-pseudobounded
topological group which is not pseudobounded. The Sorgenfrey line is not even $\omega$-pseudobounded.

\subsection{Automatic continuity of operations in semitopological groups}\label{subsec:PrelimAC}

It turns out that if the space
of a semitopological (resp. paratopological)
group satisfies suitable conditions (sometimes with some conditions imposed on the group), then the
multiplication (resp. inversion) in the group is continuous, i.e., the group is paratopological
(resp. topological). 
Finding such conditions is one of main 
branches of the theory of paratopological groups. 
 It turns out that automatic continuity essentially depends on  compact-like properties and
separation axioms of the space of a semitopological group.
 An interested reader can find known results and references on this subject
in Section 5.1 of~\cite{Rav3} and Section 3 of the survey~\cite{Tka5}
(both for semitopological and paratopological groups),
in Introduction of~\cite{AlaSan} and in the present paper (for paratopological groups).
We briefly recall the history of the topic.
In 1936 Montgomery~\cite{Mon} showed that every completely metrizable paratopological group is a
topological group. In 1953 Wallace
~\cite{Wal} asked whether every regular locally compact semitopological
group a topological group. In 1957 Ellis \cite{Ell1}, \cite{Ell2} obtained a positive
answer to the Wallace question (later the second author of the present paper showed that the regularity  condition can be relaxed,
see Proposition 5.5 in~\cite{Rav3} or its counterpart in English in~\cite{Rav6}). In 1960 \.Zelazko used Montgomery's
result and showed that each completely metrizable semitopological
group is a topological group.
Since  locally compact spaces and
completely metrizable spaces are \v{C}ech-complete (recall that
\v{C}ech-complete spaces are $G_\delta$-subspaces of compact Hausdorff spaces),
this motivated Pfister ~\cite{Pfi} in 1985 to ask whether each
\v{C}ech-complete semitopological group is a topological group. In
1996 Bouziad~\cite{Bou} and Reznichenko~\cite{Rez2}, as far as the authors know,  independently
answered affirmatively  the Pfister question. To do this, it
was sufficient to show that each \v{C}ech-complete
semitopological group is a paratopological group since earlier, Brand~\cite{Bra}
had proved that every \v{C}ech-complete paratopological group is a topological group.
Brand's proof was later improved and simplified in~\cite{Pfi}.
If $G$ is a paratopological group which is a $T_1$ space and $G\times G$ is countably compact (in
particular, if $G$ is sequentially compact) then $G$ is a topological group, see~\cite{RavRez}.
On the other hand, in the present paper we show that
we cannot weaken $T_1$ to $T_0$ here, constructing in Example~\ref{Zomega1}
a sequentially compact $T_0$ paratopological group which is not a topological group.
Also we cannot weaken the countable compactness of $G\times G$ to that
of $G$ because under additional axiomatic assumptions there exists
a countably compact paratopological group which is not a topological group,
see Example~\ref{CCNotTG}.
Also by Example~\ref{NonStandardArrowedCircle}, there exists a Hausdorff second-countable
sequentially pracompact paratopological group $G$ (so by ~\cite[Proposition 2.1]{GutRav}
each its power is sequentially pracompact) such that $G$ is not a topological group.
Remark that examples of compact-like paratopological groups that fail to be topological groups can be constructed by means of so-called cone
topologies, see Section~\ref{sec:cone-top}.
On the other hand, by Proposition~\ref{PseTG}, each
feebly compact quasiregular paratopological group is a topological group.
In particular, each pseudocompact paratopological group is a topological group,
see also~\cite[Theorem~1.7]{ArhRez} and ~\cite[Theorem~2.1]{ArhCho}.
\section{The Ra\u\i kov completion of a paratopological group}\label{s:Raikoff}
In this section we describe a canonical construction of extension of a $T_0$ paratopological group, which coincides with the Ra\u\i kov completion if the paratopological group  is a topological group.
Let us recall that a Hausdorff topological group $G$ is {\it Ra\u\i kov-complete\/} if  it is complete in the two-sided uniformity. This uniformity is generated by the family of the entourages $$\{(x,y)\in G\times G:x\in Uy\cap yU\},$$ where $U$ runs over neighborhoods of the identity in $G$. It is known that each Hausdorff topological group $G$ is a dense subgroup of a Ra\u\i kov complete topological group $\breve G$, which is unique up to a topological isomorphism. This unique topological group $\breve G$ is called {\em the Ra\u\i kov completion} of $G$, see \cite[\S3.6]{ArhTka} for more details.
In this section we extend the construction of Ra\u\i kov completion to all $T_0$ paratopological groups.
A paratopological group $(G,\tau)$ is defined to be {\em $\sharp$-complete} if its group coreflexion $(G,\tau^\sharp)$ is Hausdorff and Ra\u\i kov complete. Let us recall that a neighborhood base of the topology $\tau^\sharp$ at the identity $e$ of $G$ consists of the intersections $U\cap U^{-1}$ where $e\in U\in\tau$. The Hausdorff property of the topology $\tau^\sharp$ implies that the topology $\tau$ is $T_0$. Therefore, each $\sharp$-complete paratopological group satisfies the separation axiom $T_0$.
In this section we shall construct an embedding of any $T_0$ paratopological group $G$ into a $\sharp$-complete paratopological group $\breve G$, called the Ra\u\i kov completion of the paratopological group $G$.
The construction of the Ra\u\i kov completion exploits a construction of extension of a $\tau$-regular semigroup topology $\sigma$ on a subgroup $H$ of a semitopological group $(G,\tau)$ to a semigroup topology $\ol{\sigma}^\tau$ on $G$. This construction will be described and analyzed in Subsection~\ref{s:ext}. Its applications to Ra\u\i kov completions will be given in Subsection~\ref{s:Raikov}.
We shall need the following property of semiregularizations of paratopological groups.
\begin{lemma}\label{RegLemma} Let $(G,\tau)$ be a semitopological group, $H$ be a dense subgroup of $(G,\tau)$ and $\sigma$ be a regular semigroup topology on $H$ such that $\sigma\subseteq \tau{\restriction}H$.
Then $\sigma\subseteq \tau_{sr}{\restriction}H$. In
particular, if $\tau{\restriction}H=\sigma$ then $\tau_{sr}{\restriction}H=\sigma$.
\end{lemma}
\begin{proof}
Let $U\subseteq H$
be an arbitrary neighborhood of the identity in the topology $\sigma$.
By the regularity of $\sigma$, there exists an open neighborhood $V\in\sigma_e$ such that $\overline{V}^{\sigma}\subseteq U$.
Since $\sigma\subseteq \tau{\restriction}H$, there exists an open neighborhood $W\in\tau_e$
such that $W\cap H\subseteq V$. Taking into account that the subgroup $H$ is dense in $(G,\tau)$, we see that $\overline{W}^{\tau}=\overline{W\cap H}^{\tau}\subseteq \overline{V}^{\tau}$.
Then $\overline{W}^{\tau}\cap H\subseteq\overline{V}^{\tau}\cap H\subseteq\overline{V}^{\sigma}\subseteq U$.
\end{proof}
\subsection{Extending topologies from subgroups to groups.}\label{s:ext}
Let $H$ be a subgroup of a semitopological group $(G,\tau)$. For any semigroup topology $\sigma$ on $H$, consider the topology $\ol{\sigma}^\tau$ on $G$ consisting of the subsets $U\subseteq G$ such that for any $x\in U$ there exists a neighborhood $V\in\sigma$ of the identity such that $x\ol{V}^\tau\subseteq U$.
\begin{lemma}\label{l:neighb} For every neighborhood $U\in\sigma$ of $e$ and any $x\in G$ the set $x\ol{U}^\tau$ is a neighborhood of $x$ in the topology $\ol{\sigma}^\tau$.
\end{lemma}
\begin{proof} Since $\sigma$ is a semigroup topology, we can choose a sequence  $(U_n)_{n\in\w}$ of neighborhoods of the identity in $(H,\sigma)$  such that $U_0U_0\subseteq U$ and $U_{n}U_n\subseteq U_{n-1}$ for every $n\in\mathbb N$.  The last condition implies that
$$U_0\cdots U_n\subseteq U_0\cdots U_{n-1}U_{n-1} \subseteq U_0\cdots U_{n-2}U_{n-2}\subseteq\dots\subset
U_0U_0\subseteq U$$ for each $n\in\omega$. Taking into account that $(G,\tau)$ is a semitopological group, we conclude that
$$\overline{U_0}\cdots\overline{U_n}\subseteq \overline{U_0\cdots U_n}\subseteq \overline{U}$$ where the closure is taken in the topology $\tau$. Then the union  $W:=\bigcup_{n\in\w}\ol{U_0}\cdots\ol{U_n}$
 is a subset of $\ol{U}$ and $xW\subseteq x\ol{U}$.
We claim that the set $xW$ belongs to the topology $\ol{\sigma}^\tau$. Indeed, for any $y\in xW$, we can find $m\in\w$ such that $y\in x\ol{U_0}\cdots \ol{U_m}$ and observe that $y\ol{U_{m+1}}\subseteq x\ol{U_0}\cdots\ol{U_m}\,\ol{U_{m+1}}\subseteq xW$. Therefore, $x\ol{U}\supseteq xW\in\ol{\sigma}^\tau$ is a neighborhood of $x$ in the topology $\ol{\sigma}^\tau$.
\end{proof}
We shall prove that $\bar{\sigma}^\tau$ is a semigroup topology on $G$ if the topology $\sigma$ is {\em $\tau$-balanced} in the sense that for any neighborhood $U\in\sigma_e$ and any $x\in G$ there exists a neighborhood $V\in\sigma_e$ such that $xVx^{-1}\subseteq \overline{U}^\tau$.
A topology $\sigma$ on $H$ is defined to be {\em $\tau$-regular} if for any neighborhood $U\in\sigma_{e}$, there exists a neighborhood $V\in\sigma_e$ such that $H\cap \ol{V}^\tau\subseteq U$. Observe that a regular topology $\sigma$ on $H$ is $\tau$-regular if $\sigma\subseteq\tau{\restriction}H$.
\begin{lemma}\label{l:ext-prop} Let $H$ be a subgroup of a semitopological group $(G,\tau)$ and $\sigma$ be a $\tau$-balanced semigroup topology on $H$. Then
\begin{enumerate}
\item $\ol{\sigma}^\tau$ is a semigroup topology on $G$ such that $\ol{\sigma}^\tau{\restriction}H\subseteq \sigma$;
\item If $\sigma\subseteq\tau{\restriction}H$ and the set $H$ is dense in $(G,\tau)$, then $\ol{\sigma}^\tau\subseteq \tau$;
\item If $\sigma\subseteq\tau{\restriction}H$, then $\sigma_{sr}\subseteq\ol{\sigma}^\tau{\restriction}H$;
\item If $\sigma\supseteq\tau{\restriction}H$ and $\tau$ is a semigroup topology, then $\tau_{sr}\subseteq\ol{\sigma}^\tau$;
\item If the topology $\sigma$ is $\tau$-regular, then $\sigma=\ol{\sigma}^\tau{\restriction}H$;
\item If $\sigma\subseteq\tau{\restriction}H$ and the space $(H,\sigma)$ is regular, then $\sigma=\sigma_{sr}=\ol{\sigma}^\tau{\restriction}H$.
\item If $\sigma\subseteq\tau{\restriction}H$, the space $(H,\sigma)$ is regular, and $H$ is dense in $(G,\tau)$, then $\sigma=\sigma_{sr}=\ol{\sigma}^\tau{\restriction}H=
\left(\overline{\sigma}^\tau\right)_{sr}{\restriction}H$.
\end{enumerate}
\end{lemma}
\begin{proof} 1. The definition of the topology $\ol{\sigma}^\tau$ ensures that $\ol{\sigma}^\tau{\restriction}H\subseteq \sigma$. Next, we show that $\ol{\sigma}^\tau$ is a semigroup topology. Given any points $x,y\in G$ and a neighborhood $O_{xy}\in\bar{\sigma}^\tau$ of $xy$, use Lemma~\ref{l:neighb} to find a neighborhood $U\in\sigma$ of $e$ such that $xy\ol{U}\subseteq O_{xy}$, where the closure is taken in the topology $\tau$. Since $\sigma$ is a semigroup topology, there exists a neighborhood $V\in\sigma$ of $e$ such that $VV\subseteq U$. By the $\tau$-balanced property of the topology $\sigma$, there exists a neighborhood $W\in\sigma_e$ such that $W\subseteq V$ and $y^{-1}Wy\subseteq\ol{V}$. Since $(G,\tau)$ is a semitopological group, $y^{-1}\ol{W}y\subseteq\ol{V}$. By Lemma~\ref{l:neighb}, the sets $x\ol{W}$ and $y\ol{W}$ are neighborhoods of the points $x,y$, respectively in the topology $\ol{\sigma}^\tau$. Finally, observe that
$$x\ol{W}y\ol{W}=xy(y^{-1}\ol{W}y)\ol{W}\subseteq xy\ol{V}\cdot\ol{W}\subseteq xy\ol{VW}\subseteq xy\ol{VV}\subseteq xy\ol{U}\subseteq O_{xy},$$
witnessing that $\ol{\sigma}^\tau$ is a semigroup topology on $G$.
\smallskip
2. Assume that $\sigma\subseteq\tau{\restriction}H$ and $H$ is dense in $(G,\tau)$. To show that $\ol{\sigma}^\tau\subseteq\tau$, fix any neighborhood $W\in\ol{\sigma}^\tau_e$ of the identity and using the definition of the topology $\overline{\sigma}^\tau$, find a neighborhood $U\in\sigma_e$ such that $\ol{U}^\tau\subseteq W$. Since $\sigma\subseteq\tau{\restriction}H$, there exists a neighborhood $V\in\tau_e$ such that $V\cap H\subseteq U$. The density of $H$ in $(G,\tau)$ implies that $$V\subseteq\overline{V}^\tau=
\overline{V\cap H}^\tau\subseteq\overline{U}^\tau\subseteq W,$$witnessing that $\ol{\sigma}^\tau\subseteq\tau$.
\smallskip
3. Let $U_{sr}$ be an arbitrary neighborhood of the identity in the semiregularization $\sigma_{sr}$ of the topology $\sigma$. By Lemma~\ref{l:semiregularization}, the topology $\sigma_{sr}$ is regular. Consequently,
there exists an open neighborhood $V_{sr}\in\sigma_{sr}$ of the identity such that $\overline{V_{sr}}^{\sigma_{sr}}\subseteq U_{sr}$.
By the definition of the topology $\sigma_{sr}$, there exists a regular open set $V\in\sigma_e$  such that $V\subseteq V_{sr}$. 
By Lemma~\ref{l:neighb}, the set $\ol{V}^\tau$ is a neighborhood of the identity in the topology $\ol{\sigma}^\tau$. Since $\sigma\subseteq\tau{\restriction}H$,
$$\overline{V}^\tau\cap H\subseteq \overline{V_{sr}}^\tau\cap H=\ol{V_{sr}}^{\tau{\restriction}H}\subseteq \overline{V_{sr}}^{\sigma}\subseteq \overline{V_{sr}}^{\sigma_{sr}}\subseteq U_{sr}.$$
4. Assuming that $\tau{\restriction}H\subseteq\sigma$ and $\tau$ is a semigroup topology, we shall show that $\tau_{sr}\subseteq \overline\sigma^\tau$. By Lemma~\ref{l:semiregularization}, the topology $\tau_{sr}$ is regular. Given any neighborhood
$U_{sr}\in\tau_{sr}$ of the identity, use the regularity of the topology $\tau_{sr}$ and find a neighborhood $V_{sr}\in\tau_{sr}$ of the identity
such that $\overline{V_{sr}}^{\tau_{sr}}\subseteq U_{sr}$.
By definition of the topology $\tau_{sr}$, there exists a regular open neighborhood $V\in\tau$ such that $V\subseteq V_{sr}$. 
Since $\tau{\restriction}H\subseteq\sigma$, the set $V\cap H$ is a neighborhood of the identity in $(H,\sigma)$. By Lemma~\ref{l:neighb}, the closure $\ol{V\cap H}^\tau$ contains a neighborhood $W\in\ol{\sigma}^\tau$ of the identity. Then $W\subseteq \overline{V}^\tau\subseteq \overline{V_{sr}}^\tau\subseteq\overline{V_{sr}}^{\tau_{sr}}\subseteq U_{sr}$.
\smallskip
5. The statement follows from Lemma~\ref{l:neighb}.
\smallskip
6. If $\sigma\subseteq\tau{\restriction}H$ and the space $(H,\sigma)$ is regular, then the topology $\sigma$ is $\tau$-regular and then $\sigma=\bar{\sigma}^\tau{\restriction}H$ by the preceding statement. The regularity of the topology $\sigma$ ensures that $\sigma=\sigma_{sr}$.
\smallskip
7. Finally assume that $\sigma\subseteq\tau{\restriction}H$, the space $(H,\sigma)$ is regular, and $H$ is dense in $(G,\tau)$. By the preceding statement,  $\sigma=\sigma_{sr}=\ol{\sigma}^\tau{\restriction}H$. By  the statement (2), the density of $H$ in $(G,\tau)$ implies the density of $H$ in $(G,\ol{\sigma}^\tau)$. Now Lemma~\ref{RegLemma} implies that $\left(\overline{\sigma}^\tau\right)_{sr}{\restriction}H=\sigma$.
\end{proof}

Now we establish sufficient conditions of the $\tau$-balancedness and $\tau$-regularity of the topology $\sigma$. 

\begin{lemma}\label{l:tau-normal} Let $H$ be a subgroup of a paratopological  group $(G,\tau)$ and $\sigma$ be a semigroup topology on $H$ such that $\sigma\subseteq\tau^\sharp{\restriction}H$ and $H$ is dense in $(G,\tau^\sharp)$. Then the topology $\sigma$ is $\tau$-balanced.
\end{lemma}
\begin{proof} 
The inclusion $\sigma\subseteq\tau^\sharp{\restriction}H$ and the invariance of the topology $\tau^\sharp$ under the inversion imply that $\sigma^\sharp\subseteq\tau^{\sharp}{\restriction}H$.
Given any neighborhood $U\in\sigma_e$ and  point $x\in G$, we need to find a neighborhood $V\in\sigma_e$ such that $xVx^{-1}\subseteq\ol{U}^\tau$.
%
Since $\sigma$ is a semigroup topology, there exists a neighborhood $W\in\sigma_e$  such that $W^3\subseteq U$.  The set $W\cap W^{-1}$ belongs to the topology $\sigma^{\sharp}\subseteq\tau^\sharp{\restriction}H$. Taking into account that
the subgroup $H$ is dense in $(G,\tau^\sharp)$,
we see that the $\tau^\sharp$-closure $\overline{W\cap W^{-1}}^{\tau^\sharp}$ of the set $W\cap W^{-1}$ has nonempty interior in $(G,\tau^\sharp)$ and hence has a common point $y$ with the dense set $xH$.  It follows that $y\in\overline{W\cap W^{-1}}^{\tau^\sharp}\subseteq \overline{W}^\tau$ and $y^{-1}\in \big(\overline{W\cap W^{-1}}^{\tau^\sharp}\big)^{-1}=\overline{W\cap W^{-1}}^{\tau^\sharp}\subseteq\overline{W}^\tau$.
Since $(H,\sigma)$ is a semitopological group and $x^{-1}y\in H$, the set $V:=x^{-1}yWy^{-1}x$ is
a neighborhood of the identity in $(G,\sigma)$. Moreover,
$$xVx^{-1}=xx^{-1}yWy^{-1}xx^{-1}=yWy^{-1}\subseteq \overline{W}^\tau W\overline{W}^\tau\subseteq\overline{WWW}^\tau\subseteq\overline{U}^\tau.$$
\end{proof}
\begin{lemma}\label{l:tau-regular} Let $H$ be a subgroup of a semitopological group $(G,\tau)$. A semigroup topology $\sigma$ on $H$ is $\tau$-regular if $\sigma^{-1}\subseteq\tau{\restriction}H$, where $\sigma^{-1}=\{U^{-1}:U\in\sigma\}$.
\end{lemma}
\begin{proof} To prove that $\sigma$ is $\tau$-regular, fix any neighborhood $U\in\sigma$ of the identity. Since $\sigma$ is a semigroup topology, there exists a neighborhood $V\in\tau_e$ such that $VV\subseteq U$. Taking into account that $\sigma^{-1}\subseteq\tau{\restriction}H$, conclude that $$\ol{V}^{\tau}\cap H=\ol{V}^{\tau{\restriction}H}\subseteq\ol{V}^{\sigma^{-1}}=\bigcap_{e\in W\in\sigma^{-1}}VW^{-1}= \bigcap_{e\in W\in\sigma}VW\subseteq VV\subseteq U,$$
witnessing that the topology $\sigma$ is $\tau$-regular.
\end{proof}
\subsection{The Ra\u\i kov completion of a paratopological group}\label{s:Raikov}
Given a $T_0$ paratopological group $(G,\tau)$,
consider its group coreflexion $G^\sharp=(G,\tau^\sharp)$, which is a Hausdorff topological group.
Let $\breve G^\sharp=(\breve G,\breve \tau^\sharp)$ be the Ra\u\i kov completion
of the topological group $G^\sharp$. Observe that $\tau\subseteq\tau^\sharp=\breve\tau^\sharp{\restriction}G$. By Lemmas~\ref{l:tau-normal} and \ref{l:tau-regular}, the topology $\tau$ is $\breve\tau^\sharp$-balanced and $\breve\tau^\sharp$-regular. Consequently, the semigroup topology $\ol{\tau}^{\breve\tau^\sharp}$ on $\breve G$ is well-defined. To simplify notations, we shall denote the topology $\ol{\tau}^{\breve\tau^\sharp}$ by $\breve\tau$. The paratopological group $(\breve G,\breve\tau)$ will be called the {\em Ra\u\i kov completion} of the paratopological group $(G,\tau)$ and will be denoted by $\breve G$. By Lemma~\ref{l:neighb}, a neighborhood base at the identity in $\breve G$ consists of $\breve\tau^\sharp$-closures of neighborhoods of the identity in the paratopological group $(G,\tau)$.

Lemma~\ref{l:ext-prop}  implies the following properties of the Ra\u\i kov completions of paratopological groups. 
\begin{theorem}\label{t:Raikov} For a $T_0$ paratopological group $(G,\tau)$ its Ra\u\i kov completion $(\breve G,\breve\tau)$ has the following properties:
\begin{enumerate}
\item $\tau_{sr}\subseteq \tau=\breve\tau{\restriction}G$;
\item $\breve\tau\subseteq \breve\tau^\sharp=(\breve\tau)^\sharp$;
\item the topology $\breve\tau$ satisfies the separation axiom $T_0$;
\item If the space $(G,\tau)$ is regular, then $\tau=\breve\tau{\restriction}G=\breve\tau_{sr}{\restriction}G$, $\breve\tau^\sharp=(\breve\tau)^\sharp=(\breve\tau_{sr})^\sharp$, and the topologies $\breve\tau_{sr}\subseteq \breve\tau$ are Hausdorff;
\item If $(G,\tau)$ is a topological group, then $\breve\tau=\breve\tau^\sharp$.
\end{enumerate}
\end{theorem}
\begin{proof}
1. Lemma~\ref{l:tau-regular} implies that the topology $\tau$ is $\breve\tau^\sharp$-regular. By Lemma~\ref{l:ext-prop}(5), the $\breve\tau^\sharp$-regularity of the topology $\tau$ implies the equality $\tau=\breve\tau{\restriction}G$. The inclusion $\tau_{sr}\subseteq\tau$ follows from the definition of the semiregularization of the topology $\tau$.
\smallskip
2. The inclusion $\breve \tau\subseteq\breve\tau^\sharp$ follows from    Lemma~\ref{l:ext-prop}(2). Applying to this inclusion the operation of group coreflexion, we obtain $(\breve\tau)^\sharp\subseteq (\breve\tau^\sharp)^\sharp=\breve\tau^\sharp$. To show that $\breve\tau^\sharp\subseteq(\breve\tau)^\sharp$, take any neighborhood $U\in\breve\tau^\sharp$ of the identity. Using the regularity of the group topology $\breve\tau^\sharp$, find a neighborhood $W\in\breve\tau^\sharp$ of the identity such that $\ol{W}\subseteq U$ where the closure of $W$ is taken in the topology $\breve\tau^\sharp$. By  the definition of the topology  $\tau^\sharp=\breve\tau^\sharp{\restriction}G$, there exists a neighborhood $V\in\tau$ such that $V\cap V^{-1}\subseteq W\cap G$. Since $\tau=\breve\tau{\restriction}G$, there exists a set $T\in\breve\tau$ such that $T\cap G=V$. Since $T\cap T^{-1}\in(\breve\tau)^\sharp\subseteq\breve\tau^\sharp$, the intersection $T\cap T^{-1}\cap G$ is dense in $T\cap T^{-1}$. Consequently,
$$T\cap T^{-1}\subseteq\ol{T\cap T^{-1}\cap G}=\ol{V\cap V^{-1}}\subseteq \ol{W}\subseteq U$$
and hence $U$ is the neighborhood of the identity in the topology $(\breve\tau)^\sharp$.
\smallskip
3. To see that the topology $\breve\tau$ satisfies the separation axiom $T_0$, take any distinct points $x,y\in\breve G$. Using the $T_1$ property of the group topology  $\breve\tau^\sharp=(\breve\tau)^\sharp$, find a neighborhood $U\in\breve\tau$ of the identity such that $x^{-1}y\notin U\cap U^{-1}$. If $x^{-1}y\notin U$, then $xU\in\breve\tau$ is a neighborhood of $x$ that does not contain $y$. If $x^{-1}y\in U$, then $x^{-1}y\notin U^{-1}$ and then $yU\in\breve\tau$ is a neighborhood of $y$ that does not contain $x$.
\smallskip
4. If the space $(G,\tau)$ is regular, then  $\tau=\breve\tau{\restriction}G=\breve\tau_{sr}{\restriction}G$ by Lemma~\ref{l:ext-prop}(7). The inclusion $\breve\tau_{sr}\subseteq \breve\tau\subseteq\breve\tau^\sharp$ established in the statement (2) implies   $(\breve\tau_{sr})^\sharp\subseteq (\breve\tau)^\sharp\subseteq(\breve\tau^\sharp)^{\sharp}=\breve\tau^\sharp$. To see that $\breve\tau^\sharp\subseteq (\breve\tau_{sr})^\sharp$, fix any neighborhood $U\in\breve\tau^\sharp$ of the identity. By the regularity of the group topology $\breve\tau^\sharp$, there exists a neighborhood $W\in\breve\tau^\sharp$ of the identity such that $\overline{W}\subseteq U$, where the closure in taken in the topological group $\breve G^\sharp$. Since $W\cap G\in\tau^\sharp$, there exists a neighborhood $V\in\tau$ of the identity such that $V\cap V^{-1}\subseteq W\cap G\subseteq W$. Since $\tau=\breve\tau_{sr}{\restriction}G$, there exists a neighborhood $B\in\breve\tau_{sr}$ of $e$ such that $B\cap G\subseteq V$. Then $B\cap B^{-1}\cap G\subseteq V\cap V^{-1}\subseteq W$. Since the set $B\cap B^{-1}$ is open in the topological group $\breve G^\sharp$,
$$B\cap B^{-1}\subseteq \overline{B\cap B^{-1}}=\overline{B\cap B^{-1}\cap G}\subseteq\overline{W}\subseteq U.$$Taking into account that $U\supseteq B\cap B^{-1}\in(\breve\tau_{sr})^\sharp$, we conclude that $U$ is a neighborhood of the identity in the topology $(\breve\tau_{sr})^\sharp$ and hence   $\breve\tau^\sharp\subseteq (\breve\tau_{sr})^\sharp$. By analogy with the proof of the $T_0$ property of the topology $\breve\tau$, we can prove that the equality $\breve\tau^\sharp=(\breve\tau_{sr})^\sharp$ implies the $T_0$ separation axiom for the topology $\breve\tau_{sr}$. Being regular, the $T_0$ topology $\breve\tau_{sr}$ is Hausdorff and so is the topology $\breve\tau\supseteq\breve\tau_{sr}$.
\smallskip
5. If $(G,\tau)$ is a topological group, then $\tau=\breve\tau^{\sharp}{\restriction}G$ and $\breve\tau=\breve\tau^\sharp$ by Lemma~\ref{l:ext-prop}(2,4).
\end{proof}
Theorem~\ref{t:Raikov}(1--3) implies the following corollary.
\begin{corollary} Every $T_0$ paratopological group $G$ is a dense subgroup of the  $\sharp$-complete paratopological group $\breve G$.
\end{corollary}
Theorem~\ref{t:Raikov}(4) and Lemma~\ref{l:semiregularization} imply the following corollary.
\begin{corollary} Every $T_3$ paratopological group $G$ is a dense subgroup of the regular $\sharp$-complete paratopological group $\breve G_{sr}$, which coincides with the semiregularization of the Ra\u\i kov completion $\breve G$ of $G$.
\end{corollary}

The following simple example shows that separation axioms $T_1$ and $T_2$ are not inherited by the Ra\u\i kov completions of paratopological groups.

\begin{example} There exists a countable Abelian first-countable Hausdorff paratopological group $G$ whose Ra\u\i kov completion $\breve G$ does not satisfy the separation axiom $T_1$. Moreover,  the closure of the singleton $\{0\}$ in $\breve G$ is not compact.
\end{example}

\begin{proof} Take any countable dense subgroup $G$ in the additive group $\IR$ of real numbers. Take any non-zero real number $c$ such that $G\cap(\IZ\cdot c)=\{0\}$. For every $n\in\w$ consider the set $U_n=\bigcup_{k\in\w}\{g\in G:|g-kc|<2^{-n}\}$. Observe that $U_{n+1}+U_{n+1}\subseteq U_n$ for every $n\in\w$, which implies that the topology 
$$\tau=\big\{U\subseteq G:\forall x\in U\;\exists n\in\w\quad (x+U_n\subseteq U)\}$$ turns $G$ into a first-countable paratopological group.  The paratopological group $(G,\tau)$ is Hausdorff since for any distinct elements $x,y\in G$ there exists $n\in\w$ such that $x-y\notin \bigcup_{k\in\IZ}\{g\in G:|g-kc|<2^{-n}\}$, which implies that $x+U_{n+1}$ is disjoint with $y+U_{n+1}$. The topology $\tau^\sharp$ of the group coreflexion  of the paratopological group $(G,\tau)$ coincides with the Euclidean topology on $G$ and the Ra\u\i kov completion of the topological group $(G,\tau^\sharp)$ can be identified with the real line. The definition of the topology $\tau$ ensures that for any $n\in\w$ the closure of the set $U_n$ in the real line contains the set $\{kc:k\ge 0\}$. Consequently, the closure of $\{0\}$ in the Ra\u\i kov completion $\breve G$ of $(G,\tau)$ contains the set $\{kc:k\le 0\}$ and hence is not compact.
\end{proof}

\begin{remark} After writing down the above construction of the Ra\u\i kov completion of a $T_0$ paratopological group, we have discovered that essentially the same completion was elaborated by Mar\'\i n and Romaguera \cite{MarRom1}  who observed that the bicompletion of the two-sided quasiuniformity of a paratopological group  has a natural structure of a paratopological group that  contains the group as a 2-dense subgroup. In contrast, our construction is more direct and does not involve quasiuniformities.
\end{remark}

\section{When is a feebly compact paratopological group a topological group?}\label{PPGTAT}
The following lemma proved by Tkachenko \cite[Theorem 3.1]{Tka7} allows us to reduce topological questions related to
semitopological group to the case of $T_0$-spaces.
\begin{lemma}\label{l:Tka-T0} Let $G=(G,\tau)$ be a semitopological group and $H:=\bigcap_{U\in\tau_e}(U\cap U^{-1})$. Then $H$ is a normal subgroup of $G$, the quotient space $T_0G:=G/H$ satisfies the separation axiom $T_0$ and each open set $U\subseteq G$ is equal to $\pi^{-1}(\pi(U))$ where $\pi:G\to G/H$ is the quotient homomorphism.
\end{lemma}
The quotient group  $T_0G:=G/H$ is called the {\em $T_0$-reflexion} of $G$. If $G$ is a (para)topological group, then so is its $T_0$-reflexion $T_0G$, see \cite{Tka7}.

\begin{proposition}\label{TCCTG} Each totally countably compact paratopological group is a topological group.
\end{proposition}
\begin{proof}
Let $G=(G,\tau)$ be such a group. Put $B=\bigcap\tau_e$. The continuity of the multiplication implies that $B$ is a subsemigroup of the group $G$. Put $B'=\ol{\{e\}}=\{x\in G:\forall U\in\tau_e\;\;(e\in xU)\}=\{x\in G:e\in xB\}=B^{-1}$.
Hence $B'$ is a closed subsemigroup of the group $G$. The total countable compactness of $G$ implies that $B'=\overline{\{e\}}$ compact. Being a compact
topological semigroup, $B'$ contains a minimal closed right ideal $I\subset
B'$. Since $I$ is closed in $G$, for any element $x\in I$, the set $x^2I\subseteq B'$ is closed in $G$ and in $B'$. By the minimality of $I$, $x^2I=I\ni x$.
Hence $x^{-1}\in I$ and $e\in I$. Consequently $I=B'$ and $xB'=B'$ for
each element $x\in B'$. Therefore $B'$ is a group and $B'=B$.
Since $B=\bigcap\tau_e$, we see that $g^{-1}Bg\subseteq B$ for each $g\in
G$. Hence $B$ is a normal subgroup of the group $G$.
Since $B$ has the antidiscrete topology, $B$ is a topological group.
Since the normal subgroup $B$ is closed in $G$, the quotient group $G/B$ is a $T_1$-space. Let $\pi:G\to G/B$ be the quotient homomorphism.
Given any sequence $(y_n)_{n\in\omega}$ in the quotient group $G/B$, choose a sequence $(x_n)_{n\in\omega}$ in the group $G$ such that $\pi(x_n)=y_n$
for each $n\in\omega$. Since $G$ is totally countable compact,  there is a subsequence $A$ of
$\{x_n:n\in\omega\}$ such that the closure $\ol A$ is compact. Since the set $\ol A$ is closed, $\ol A=\ol AB$.
Hence the closed compact set $\pi(\ol A)$ contains a subsequence of the sequence $\{y_n:n\in\omega\}$.
Therefore $G/B$ is a $T_1$ totally countably compact paratopological group. By  \cite[Cor 2.3]{AlaSan}, $G/B$ is a topological group. Hence both $B$ and $G/B$ are topological groups.
Then by \cite[1.3]{Rav4} or \cite[Proposition 5.3]{Rav3}, $G$ is a topological group too.
\end{proof}
Proposition~\ref{TCCTG} generalizes Corollary 2.3 in \cite{AlaSan} and Lemma 5.4 in \cite{Rav3}.
\begin{lemma}\label{SatQua} \cite[Proposition 1.7]{Rav2} Each saturated paratopological group is
quasiregular.
\end{lemma}
Let us recall \cite{BokGur} that a paratopological group $G$ is {\it topologically
periodic} if for each $x\in G$ and a neighborhood $U\subseteq G$ of
the identity there is a number $n\ge 1$ such that $x^n\in U$.
Reznichenko proved the following lemma, announced in~\cite{RavRez}.

\begin{lemma}\label{RegTG}  A paratopological group $(G,\tau)$ is a topological group if and only if $(G,\tau)$ is quasiregular and its semiregularization $(G,\tau_{sr})$ is a topological group.
\end{lemma}
\begin{proof} The ``only if'' part follows trivially from the regularity of any topological group. To prove the ``if'' part, assume that $(G,\tau)$ is quasiregular and its semiregularization $(G,\tau_{sr})$ is a topological group. To show that $G$ is a topological group, fix any neighborhood $N\in\tau$ of $e$. Since $(G,\tau)$ is a paratopological group, there exists a neighborhood $U\in\tau$ of $e$ such that $UU\subseteq N$. Consider the regular open set $W:=\int\ol U\in\tau_{sr}$ and observe that the open set $V=\bigcup\{T\in\tau: \ol{T}\subseteq U\}$ is dense in $U$ (by the quasiregularity of the topology $\tau$) and belongs to the topology $\tau_{sr}$, being the union of canonically open sets in the topology $\tau$. Consequently, $V\subseteq U\subseteq W$ and $\ol V=\ol U=\ol W$. Since
$(G,\tau_{sr})$ is a topological group, the sets $W^{-1}$ and $V^{-1}$ are $\tau_{sr}$-open. The
density of the set $V$ in $W$ implies the density of the set $V\cap V^{-1}$ in $W\cap W^{-1}$ in
the topological group $(G,\tau_{sr})$. Since $W\cap W^{-1}\ni e$ is nonempty, the open set $V\cap
V^{-1}\in\tau_{sr}\subseteq\tau$ is not empty, too. Take any point $x\in V\cap V^{-1}$ and observe
that $$N^{-1}\supseteq (UU)^{-1}\supseteq V^{-1}V^{-1}\supseteq x^{-1}(V\cap V^{-1})$$ is a neighborhood
of $e$ in the topology $\tau$.
\end{proof}
The following lemma easily follows from the definitions.
\begin{lemma}\label{RegPCompact} A space $(X,\tau)$ is feebly compact  if and only if its semiregularization $(X,\tau_{sr})$ is feebly compact.\qed
\end{lemma}
The following two lemmas were proved by Alas and Sanchis \cite{AlaSan} in case of paratopological groups satisfying the separation axiom $T_0$. But Lemma~\ref{l:Tka-T0} allows us to remove the $T_0$ assumption.
\begin{lemma}\label{TopPerSat} Each topologically periodic Baire paratopological group is saturated.
\end{lemma}
\begin{lemma}\label{2PseBai} Each $2$-pseudocompact paratopological group is a Baire space.
\end{lemma}
\begin{lemma}\label{PseudoTGPre} Each feebly compact topological group $G$ is precompact.\qed
\end{lemma}
\begin{proof} If $G$ is not precompact, then there exists a neighborhood $U\subseteq G$ of the identity and a sequence of points $\{x_n\}_{n\in\omega}\subseteq G$ such that $x_m\notin x_nU$ for any $n<m$. Take any neighborhood $V=V^{-1}\subseteq G$ of the identity with $V^4\subseteq U$ and observe that $\{x_nV\}_{n\in\omega}$ is an infinite locally finite family of open sets in $G$, witnessing that $G$ is not feebly compact.
\end{proof}
The following lemma was first proved in the Ph.D. dissertation \cite[3.1]{Rav3} of the first author.
\begin{lemma}\label{PreSat} Every precompact paratopological group $(G,\tau)$ is saturated.\qed
\end{lemma}
\begin{proof}  Given any  neighborhood $U\in\tau$ of the identity $e$, choose a neighborhood $V\in\tau$ of $e$ such that $VV\subseteq U$. By the precompactness of $(G,\tau)$, there exists a finite set $F\subseteq G$ such that $G=FV$. Since $G=V^{-1}F^{-1}$, for some $x\in F^{-1}$ the set $V^{-1}x^{-1}$ is not nowhere dense in $G$. Then the closure $\ol{V^{-1}}$ contains some nonempty open set $W\in\tau$. Now observe that $$W\subseteq \ol{V^{-1}}=\bigcap_{T\in\tau_e}V^{-1}T^{-1}\subseteq V^{-1}V^{-1}\subseteq U^{-1},$$which means that the set $U^{-1}$ has nonempty interior. Therefore, the paratopological group $(G,\tau)$ is saturated.
\end{proof}
By analogy one can prove the following lemma, see \cite[Theorem~2.5]{AlaSan}.
\begin{lemma}\label{BaiSat}
Each Baire left $\omega$-precompact paratopological group is saturated.\qed
\end{lemma}
The following lemma can be derived from \cite[Theorem~1.7]{ArhRez} and Lemma~\ref{l:Tka-T0}.
\begin{lemma}\label{ArhRezT3}
Each paratopological group that is a dense $G_\delta$-subset of a regular feebly compact space is a topological group.
\end{lemma}
\begin{lemma}\label{2PseT3}
Each regular $2$-pseudocompact paratopological group is a topological group.
\end{lemma}
\begin{proof} Let $(G,\tau)$ be such a group.
Suppose $(G,\tau)$ is not a topological group. Then there exists a neighborhood $U\in\tau$ of the unit $e$ such that
$V\not\subseteq\ol{U^{-1}}\subseteq (U^{-1})^2$ for every $V\in\tau_e$.
By induction we can build a sequence
$\{U_i\}_{i\in\omega}\in\tau$ of open neighborhoods of $e$ such that
$\ol{U_0}\subseteq U$ and $U_{i+1}^2\subseteq U_i$ for each $i\in\omega$. Then $F=:\bigcap_{i\in\w}U_i$ is a subsemigroup of $G$. Since the group $G$ is $2$-pseudocompact, there is a point $x\in G$ such that
$$x\in\bigcap_{i\in\w}\ol{(U_i\bs\ol{U^{-1}})^{-1}}\subseteq
\bigcap_{i\in\w}\ol{U_i^{-1}}\subseteq\bigcap_{i\in\w} (U_i^{-1})^2\subseteq F^{-1}.$$ Moreover,
$x\in\ol{(U_0\bs\ol{U^{-1}})^{-1}}\subseteq
\ol{(U_0\bs U^{-1})^{-1}}\subseteq \ol{U_0^{-1}\bs U}\subseteq G\bs U$.
Since $x\not\in\ol {U_0}$ there is a neighborhood $W\in\tau_e$ such that
$W\subseteq U_1$ and $Wx\cap U_0=\0$.
Let $m,n\in\omega$ and $n>m\ge 1$. If $\emptyset\not=Wx^n\cap Wx^m$ then
$\emptyset\not=Wx \cap Wx^{m+1-n}\subseteq Wx\cap U_1F\subseteq Wx\cap U_1^2\subseteq Wx\cap U_0=\emptyset$.
Therefore $\{Wx^i:i\ge 1\}$ is an infinite family of pairwise disjoint sets.
Choose a neighborhood $V\in\tau_e$ such that $V^2\subseteq W$.
Since $G$ is $2$-pseudocompact, there is
a point $y\in\bigcap_{n\ge 1}\ol{(\bigcup_{i\ge n} Vx^i)^{-1}}$.
There exist numbers $m,n\in\omega$ such that $n>m\ge 1$ and $yV\cap (Vx^m)^{-1}\not=\emptyset\ne yV\cap (Vx^n)^{-1}$. Then $y^{-1}\in V^2x^m\cap V^2x^n\subseteq Wx^m\cap Wx^n$.
This contradiction proves that $G$ is a topological group.
\end{proof}
Observe that a paratopological group is $2$-pseudocompact if it is either countably compact or it is a feebly compact topological group.
\begin{proposition}\label{2PseP}
Each $2$-pseudocompact paratopological group is feebly compact.
\end{proposition}
\begin{proof} Let $G$ be such a group.
The semiregularization  $G_{sr}$ of $G$ is a regular $2$-pseudocompact
paratopological group. By Lemma~\ref{2PseT3}, $(G,\tau_{sr})$ is a topological group. The $2$-pseudocom\-pactness of the topological group $G_{sr}$ implies the feeble compactness of $G_{sr}$. By Lemma~\ref{RegPCompact}, the paratopological group $G$ is feebly compact.
\end{proof}

Proposition~\ref{2PseP} cannot be reversed: in \cite[Theorem~2]{SanMTka} Tkachenko and Sanchis have constructed a Hausdorff feebly compact Baire paratopological group which is not $2$-pseudocompact, thus answering a question of the second author from an older version of this manuscript.
\begin{problem}\label{PCto2PC} Is each countably pracompact Baire paratopological group
$2$-pseudocompact?
\end{problem}
The answer to Problem~\ref{PCto2PC} is positive if the paratopological group is quasiregular or saturated or topologically periodic or left $\omega$-precompact.

\begin{proposition}\label{PseTG} For a  feebly compact paratopological group $G$,   the following conditions
are equivalent:
\begin{itemize}
\item[\textup{1)}] $G$ is quasiregular;
\item[\textup{2)}] $G$ is saturated;
\item[\textup{3)}] $G$ is topologically periodic and Baire;
\item[\textup{4)}] $G$ is left $\omega$-precompact and Baire;
\item[\textup{5)}] $G$ is precompact;
\item[\textup{6)}] $G$ is a topological group.
\end{itemize}
\end{proposition}
\begin{proof} It is easy to prove that Condition 6 implies all other conditions
(remark that a feebly compact topological group is $2$-pseudocompact and hence Baire by Lemma~\ref{2PseBai}).
\smallskip
$(1)\Rightarrow(6)$ By Lemma~\ref{l:semiregularization}, the semiregularization  $G_{sr}$ of the paratopological group $G$ is a regular feebly compact
paratopological group. By Lemma~\ref{ArhRezT3}, $G_{sr}$ is a topological group.
Then by Lemma~\ref{RegTG}, $G$ is a topological group, too.
The implication $(2)\Rightarrow (1)$ follows from Lemma~\ref{SatQua}.
The implications $(3,4,5)\Rightarrow (2)$ follow from Lemmas~\ref{TopPerSat}, \ref{BaiSat}, and \ref{PreSat}, respectively.
\end{proof}
The implication $(1)\Rightarrow (6)$ in Proposition~\ref{PseTG} generalizes Proposition 2 in \cite{RavRez}; $(4)\Rightarrow (6)$
generalizes Proposition 2.6 in \cite{AlaSan}; the implication $(3)\Rightarrow (6)$
answers Question C in \cite{AlaSan} and generalizes Theorem 3 in \cite{BokGur}.

A group has two basic operations: the multiplication and the inversion. The first operation is
continuous on a paratopological group, but the second can be discontinuous. But there are
continuous operations on a paratopological group besides the multiplication. For instance, a power.
Now suppose that there exists a nonempty open subset $U$ of a paratopological group $G$ such that
the inversion on the set $U$ coincides with a continuous operation. Then the continuity of the
inversion on $U$ can be used to show that the group $G$ is topological. Let us consider a trivial
application of this idea. Let $G$ be a paratopological group of bounded exponent. Then there exists
a number $n\in\mathbb N$ such that the inversion on the group $G$ coincides with the $n$-th power.
Thus the group $G$ is topological. Another applications are formulated in two next propositions.
Let us recall that a group  is {\em periodic} if each element of the group has finite order.
\begin{proposition}\label{p:2.16} Each Baire periodic paratopological group is a topological group.
\end{proposition}
\begin{proof} Let $(G,\tau)$ be such a group. Put $B'=\bigcap\{U^{-1}:U\in\tau_e\}$. Then $B'$ is closed subset of $G$.
Since $B'$ is a periodic semigroup, we see that $B'$ is a group such that $B'=(B')^{-1}\subseteq U$ for any neighborhood $U$ of the identity.
For each $n>1$ we put $G_n=\{x\in G:x^n\in B'\}$. Since the space $(G,\tau)$ is Baire, there is a
number $n>1$ such that the set $G_n$ has the nonempty interior in $(G,\tau)$.
Therefore there are a point $x\in G_n$ and a neighborhood $U\in\tau_e$ such that
$(xu)^n\in B'$ for each $u\in U$. Then $V^{-1}\subseteq B'(xV)^{n-1}x$ for each subset $V$ of $U$.
Let $U_1\in\tau$ be an arbitrary neighborhood of the identity. Then $x^n\in B'U_1$.
By the continuity of the multiplication there is a neighborhood $W\subseteq U$ of the identity such that
$(xW)^{n-1}x\subseteq B'U_1$. Then $W^{-1}\subseteq B'(xW)^{n-1}x\subseteq B'B'U_1=B'U_1$. Taking into account that $B'\subseteq U_1$, we see that $W^{-1} \subseteq B'U_1\subseteq U_1^2$.
\end{proof}
\begin{proposition}\label{T2PPTG} Any Hausdorff feebly compact periodic paratopological
group is a topological group.
\end{proposition}
\begin{proof} Let $(G,\tau)$ be such a group. By Lemma~\ref{l:semiregularization}, its semiregularization $(G,\tau_{sr})$ is a Hausdorff regular paratopological group.  By Proposition~\ref{PseTG}, $(G,\tau_{sr})$ is a
pseudocompact periodic topological group. Let $\breve G$ be the Ra\u\i kov completion of the topological group $(G,\tau_{sr})$. The feeble compactness of $(G,\tau_{sr})$ implies the precompactness of $(G,\tau_{sr})$ and the compactness of $\breve G$.
We claim that the group $\breve G$ is periodic. Suppose to the contrary that there exists a
non-periodic element $x\in\breve G$. Then for each positive integer $n$ there exists a neighborhood
$V_n\ni x$ such that $V_n^n\not\ni e$. Being pseudocompact, the subgroup $G$ of $\breve G$
intersects each nonempty $G_\delta$-set in $\breve G$. Consequently, there exists a point $y\in
G\cap \bigcap_{n\in\w} V_n$. By the periodicity of $G$, there is a positive integer $n$ such that
$y^n=e$. But this contradicts $e\notin V_n^n$. This contradiction proves the periodicity of the
group $\breve G$.
For each $n>1$ we put $\breve G_n=\{x\in\breve G:x^n=e\}$. Since the compact topological group $\breve G$ is Baire, there is a
number $n>1$ such that the set $\breve G_n$ has the nonempty interior. Since the group
$G$ is dense in $\breve G$, we see that there are a point $x\in G$ and a neighborhood $U$ of the identity of $G$
such that $(xu)^n=e$ for each $u\in U$.
Let $y$ be an arbitrary point of $G$ and $V\in\tau$ be an arbitrary neighborhood of the point $y^{-1}$.
Then $x^ny^{-1}=y^{-1}\in V$. By the continuity of the multiplication there is a neighborhood
$W\in\tau$ of the identity such that $W\subseteq U$ and $(xW)^{n-1}xy^{-1}\in V$. Let $z\in yW$ be an
arbitrary point. Then $(xy^{-1}z)^n=e$ and hence $z^{-1}=(xy^{-1}z)^{n-1}xy^{-1}\in (xW)^{n-1}xy^{-1}\subseteq V$.
Thus the inversion on the group $(G,\tau)$ is continuous and $(G,\tau)$ is a topological group.
\end{proof}
The following example shows that Proposition~\ref{T2PPTG} cannot be generalized to $T_1$ groups
(even to sequentially pracompact and periodic).
Moreover, in \cite{AlaSan}, Alas and Sanchis proved that each topologically periodic
$T_0$ paratopological group is $T_1$. Example~\ref{BanakhEx} combined with Proposition~\ref{T2PPTG}  show that there is a
$T_1$ periodic paratopological group $G$ which is not Hausdorff.
\begin{example}\label{BanakhEx}  There exists a $T_1$ periodic sequentially pracompact
paratopological group $G$ which is not a topological group.
\end{example}
\begin{proof}
For each positive integer $n$ let $C_n$ be the cyclic group $\{0,\dots,n-1\}$ endowed with the discrete
topology and the binary operation $\oplus$ such that $x\oplus y\equiv x+y \mod n$ for any  $x,y\in C_n$.
Let $G=\bigoplus\limits_{n=1}^\infty C_n$ be the direct sum of the cyclic groups $C_n$. Elements of $G$ can be thought as functions $x:\mathbb N\to \omega$ such that $x(n)<n$ for all $n\in\mathbb N$ and $x(n)=0$ for all but finitely many $n$.
Let $\mathcal F$ be the family of all non-decreasing unbounded functions from $\mathbb N$ to
$\omega$. For each $f\in\mathcal F$ put $$O_f=\{0\}\cup\{x\in G:\exists n\in\mathbb N\; (0<x(n)<f(n)\mbox{ and }\forall m>n\;x(m)=0)\}.$$ It is easy to check that the family $\{O_f:f\in\mathcal F\}$ is a base
at the zero of a $T_1$ semigroup topology $\tau$ on $G$.
For each positive integer $m\ge 2$ consider the function $\delta_m\in G$ such that $\delta_m^{-1}(1)=\{m\}$ and $\delta_m^{-1}(0)=\mathbb N\setminus\{m\}$.
Let $f,g\in\mathcal F$ be arbitrary functions and let $x\in G$ be an arbitrary element.
There exists a number $m$ such that $f(m)\ge 2$, $g(m)\ge 2$, and $x(n)=0$ for each $n\ge m$.
Then $\delta_m\in O_g$ and $x+\delta_m\in O_f$. Therefore $x\in\overline{O_f}$. Hence $\overline{O_f}=G$ for each $f\in\mathcal
F$ and the topology $\tau$ is not Hausdorff. Since any two nonempty open subsets of $(G,\tau)$ intersect, any two nonempty open subsets of each power of $G$ in the box topology intersect, too.
Therefore each power of $G$ in the box topology is feebly compact.
Put $D=\{0\}\cup \{\delta_m:m\ge 2\}$. Then $D$ witness that $G$ is sequentially pracompact,
because $D$ is homeomorphic to a convergent sequence and dense in $G$.
\end{proof}
\begin{problem}\label{prob2} Suppose a Hausdorff topologically periodic paratopological group $G$ satisfies one of the following
conditions:
\begin{itemize}
\item $G$ is countably pracompact (and left $\omega$-precompact),
\item $G$ is feebly compact (and left $\omega$-precompact).
\end{itemize}
Is $G$ a topological group?
\end{problem}
If a paratopological group $G$ satisfies one of the conditions of Problem~\ref{prob2} and is quasiregular or saturated or Baire, then $G$ is a topological group.
We recall the following question by van Douwen~\cite{vDou} from 1980,
which is, according to references from~\cite[p. 2]{HvMR-GS},
considered central in the theory of topological groups.
\begin{question} Is there an infinite countably compact group without non-trivial convergent sequences?
\end{question}
For a history of related results see introductions of~\cite{Tom3} and~\cite{HvMR-GS}. In particular,
in 1976 Hajnal and Juh\'asz in~\cite{HJ} constructed the first example of such a group
under the Continuum Hypothesis. In 2004 Garc\'ia-Ferreira, Tomita, and Watson in~\cite{G-FTW}
obtained such a group from a selective ultrafilter. In 2009 Szeptycki and Tomita in~\cite{ST}
showed that in the Random model there exists such a group, giving the first example
that does not depend on selective ultrafilters.
Finally, Hru\v s\'ak, van Mill, Ramos-Garc\'ia, and Shelah obtained a positive answer to
van Douwen's question in ZFC, see~\cite{HvMR-GS} for ``a presentable draft of the proof".
The built group is Boolean and at p.17 the authors say that they do not know how
to modify (in ZFC) their construction to get a non-torsion example and formulate
the question about an existence of such a group in ZFC.
We denote by TT a bit stronger axiomatic assumption: {\em there is an infinite torsion-free Abelian
countably compact topological group without non-trivial convergent sequences.}
The first example of such a group was
constructed by Tkachenko~\cite{Tka4}
under the Continuum Hypothesis. Later, the Continuum Hypothesis was weakened  to Martin's Axiom
for $\sigma$-centered posets by Tomita in~\cite{Tom2}, for countable posets in~\cite{KosTomWat},
and finally to the existence of $\mathfrak c$ many incomparable selective ultrafilters in~\cite{MadTom}.
For a subset $A$ of a group $G$ by $\langle A\rangle\subseteq G$ we denote the  subgroup generated by the set $A$ in $G$.
The proof of Lemma 6.4 in \cite{BanDimGut} implies the following
\begin{lemma} \label{FA} Under TT, the free Abelian group $G$ generated by the cardinal $\mathfrak c$ admits a Hausdorff group topology $\tau$
such that for each countable infinite subset $M\subseteq G$ there exists an ordinal
$\alpha \in \mathfrak c\cap \overline M^\tau$ such that $M\subseteq\langle[0,\alpha)\rangle$.\qed
\end{lemma}
\begin{example} \label{CCNotTG} Under TT there exists a functionally Hausdorff countably compact free Abelian paratopological group $(G,\sigma)$, which is not  a topological group.
\end{example}
\begin{proof}
Let $G$ be a free Abelian group generated by the cardinal $\mathfrak c$. Elements of the group $G$ can be thought as functions $x:\mathfrak c\to\mathbb Z$ with finite support $\supp(x):=\{\alpha\in\mathfrak c:x(\alpha)\ne 0\}$. Each ordinal $\alpha\in\mathfrak c$ can be identified with its characteristic function $\delta_\alpha:\mathfrak c\to\{0,1\}$ (such that $\delta_\alpha^{-1}(1)=\{\alpha\}$).
In the free Abelian group $G$, consider the subsemigroup  $$S=\{0\}\cup \big\{x\in G\setminus\{0\}:x\big(\max(\supp(x))\big)<0\big\}.$$ By  Lemma \ref{FA}, the group $G$ admits a Hausdorff group topology $\tau$
such that for each countable infinite subset $M\subseteq G$ there exists an ordinal $\alpha \in \mathfrak c\cap \overline M^\tau$ such that
$M\subseteq \langle [0,\alpha) \rangle$.  It is easy to check that the family $\{U\cap S:0\in U\in\tau\}$
is a base at the zero of a semigroup topology $\sigma$ on $G$.
Let $M$ be an arbitrary countable infinite subset of the group
$G$. By the choice of the topology $\tau$, there exists an ordinal $\alpha \in \overline{M}^\tau\cap\cc$ such that
$M\subseteq \langle[0, \alpha) \rangle$. We claim that $\alpha$ belongs to the closure of the set $M$ in the topology $\sigma$. Given any open neighborhood $U\in\tau$ of zero, we need to prove that the set $\alpha+(U\cap S)$ intersects $M$. Since $\alpha\in\overline{M}^\tau$, there exists an element $\mu\in M\cap(\alpha+U)$. On the other hand, $\mu\in M\subseteq\langle[0,\alpha)\rangle\subseteq \alpha+S$ and hence $\mu\in (\alpha+U)\cap(\alpha+S)=\alpha+(U\cap S)$.
Therefore, $(G,\sigma)$ is a countably compact
Hausdorff paratopological group. In particular, the space $(G,\sigma)$ is not discrete. Taking into account that $S\in\sigma$ and $S\cap(-S)=\{0\}\notin\sigma$, we see that $(G,\sigma)$ is not a topological group.
\end{proof}
Example~\ref{CCNotTG} negatively answers Problem 1 from \cite{Gur} under TT.
Since each countably compact left $\omega$-precompact paratopological group is a topological group \cite[p.
82]{Rav3}, we see that Example \ref{CCNotTG} implies the negative answers to Question A in \cite{AlaSan} and
to Problem 2 in \cite{Gur} under TT.

\begin{problem}\label{prob1} Is there a ZFC-example of a paratopological group $G$ such that $G$ is not a topological
group but $G$ satisfies one of the following conditions:
\begin{itemize}
\item[\textup{(i)}] $G$ is Hausdorff countably compact,
\item[\textup{(ii)}] $G$ is Hausdorff countably pracompact Baire,
\item[\textup{(iii)}] $G$ is $T_1$ countably compact?
\end{itemize}
\end{problem}
If $G$ is such a group, then $G$ is not topologically periodic, not saturated, not quasiregular,
and not left $\omega$-precompact. Moreover, $G\times G$ is not countably compact.
By Theorem 2.2 from \cite{AlaSan}, the group $(G,\sigma)$ in Example~\ref{CCNotTG}
is Baire. Therefore under TT there is a group $G$ such that $G$ satisfies all conditions listed in Problem~\ref{prob1}.

The following proposition was published in 2010 in the first arXiv version of this manuscript
and answered Question B from~\cite{AlaSan}. In 2012 S\'anchez generalized this result,
obtaining a characterization of 2-pseudocompact $T_1$ paratopological groups which
are topological groups, see ~\cite[Proposition 2.9]{San}. This characterization was
independently obtained by Xie and Lin in~\cite{LinXie2}.
So we skip the proof of the proposition, but keep its statement for references farther in the paper.
Also this and the next proposition positively answer
a special case (for $2$-pseudocompact groups) of Problem 8.1 in \cite{ArhChoKen}.
\begin{proposition}\label{CountPseudoCh} Each $2$-pseudocompact paratopological group of countable pseudocharacter
is a topological group.
\end{proposition}
In order to strengthen Proposition~\ref{CountPseudoCh} for Hausdorff $2$-pseudocompact
paratopological groups we shall use the following lemmas.
\begin{lemma}\label{CancComp}\cite[Theorem~0.5]{Rez} A compact Hausdorff semigroup with separately continouous
multiplication and two-side cancellations is a topological group.
\end{lemma}
\begin{lemma}\label{MySwellingLemma} Let $K$ be a compact nonempty subset of a Hausdorff semitopological group
$G$. Then the set $G_K=\{x\in G: xK\subseteq K\}$ is a compact Hausdorff topological group.
\end{lemma}
\begin{proof} Observe that $G_K=\bigcap_{x\in K}Kx^{-1}$ is a compact  cancellative subsemigroup of the Hausdorff semitopological group $G$.  By Lemma~\ref{CancComp}, $G_K$ is a topological group.
\end{proof}
\begin{remark} Lemma~\ref{MySwellingLemma} does not necessarily hold for $T_1$ paratopological groups,
as the following example shows. Let $G=(\mathbb Z,+)$ be the group of integers endowed with a
topology with the base $\{\{x\}\cup \{z\in\mathbb Z:z\ge y\}:x,y\in\mathbb Z\}$. It is easy to
check that $G$ is a $T_1$ paratopological group. Put $K=\{x\in\mathbb Z:x\ge 0\}$. Then $K$ is a compact
subset of the group $G$ and $G_K=K$ is not a group.
\end{remark}
\begin{lemma}\label{SquareCompNeighb}\cite[Proposition 1.8]{Rav} If $G$ is a paratopological group,
$K$ is a compact subset of $G$, $F$ is a closed subset of $G$ and $K\cap F=\emptyset$,
then there exists a neighborhood $V\subseteq G$ of the identity such that $VK\cap F=\emptyset$.
\end{lemma}
\begin{proposition}\label{CountPseudoChComp} Each $2$-pseudocompact Hausdorff paratopological group
containing a non\-empty compact $G_\delta$-set is a topological group.
\end{proposition}
\begin{proof} Let $(G,\tau)$ be such a group and $K$ be a nonempty compact $G_\delta$-set of the space $(G,\tau)$. Suppose $(G,\tau)$ is not a topological group;
then there exists a neighborhood $U\in\tau_e$ such that
$V\not\subseteq\ol{U^{-1}}\subseteq (U^{-1})^2$ for each neighborhood $V\in\tau_e$.
By induction using Lemma~\ref{SquareCompNeighb} we can build a sequence
$\{V_i\}_{i\in\w}\subseteq\tau$ of neighborhoods of $e$ such that $V_0=U$,
$V_{i+1}^2\subseteq V_i$ for each $i\in\omega$ and $\bigcap_{i\in\omega}V_iK=K$.
Put $H=\bigcap_{i\in\omega}V_i$. The choice of the sequence $(V_i)_{i\in\w}$ ensures that $H$ is a semigroup.
 Moreover, since $\overline{V_{i+1}^{-1}}\subseteq V_{i+1}^{-1}V_{i+1}^{-1}\subseteq V_{i}^{-1}$ for each
$i\in\omega$, we see that $H^{-1}$ is a closed subset of $G$. Moreover,
$HK\subseteq \bigcap_{i\in\omega}V_iK=K$, so $H\subset
G_K$. Since $G_K$ is a group (see Lemma~\ref{MySwellingLemma}), $H^{-1}\subseteq G_K$ too, so $H^{-1}$ is a
compact cancellative topological semigroup. By Lemma~\ref{CancComp},
$H^{-1}$ is a group.
We have
$V_i\not\subseteq\ol{U^{-1}}$ for each $i\in\omega$. Since the group $G$ is
$2$-pseudocompact, there is a point $x\in G$ such that
$xW^{-1}\cap (V_i\bs\ol{U^{-1}})\not=\0$ for each $W\in\tau_e$ and each $i>0$.
Then $xW^{-1}\cap (V_i\bs U^{-1})\not=\0$ and $Wx^{-1}\cap (V_i^{-1}\bs U)\not=\0$.
Therefore $x^{-1}\in\ol{V^{-1}_i}\subseteq V_{i-1}^{-1}$. Hence
$x\in H^{-1}$. But $H^{-1}=H\subseteq U\cap U^{-1}$ and $U^{-1}\cap (V_1\bs\ol{U^{-1}})=\0$. This contradiction proves
that $G$ is a topological group.
\end{proof}
The authors do not know whether counterparts of Proposition~\ref{CountPseudoChComp} hold for $T_1$
2-pse\-u\-do\-compact or countably compact paratopological groups.
But a $T_0$ sequentially compact paratopological group $G_S$ from Example~\ref{Zomega1} contains an open compact subsemigroup $S$. Nevertheless, $G_S$ is not a topological group.
The following Example~\ref{NonStandardArrowedCircle} (which is an extension of  \cite[Example 5.17]{Rav3}) shows that counterparts
of Proposition \ref{CountPseudoCh} and Implication $(4)\Rightarrow (1)$ in Proposition~\ref{PseTG}
do not hold for sequentially pracompact paratopological groups.
Also Example~\ref{NonStandardArrowedCircle} negatively answers Problem 6.2 in \cite{Tka5}
and Problem 8.1  in \cite{ArhChoKen}.
\begin{example}\label{NonStandardArrowedCircle}
There exists a Hausdorff second countable sequentially pracompact
paratopological group $G$ which is not a topological group.
\end{example}
\begin{proof} Let $\IT=\{z\in\IC:|z|=1\}$ be the unit circle endowed with the standard Euclidean topology $\tau$ and the group operation of multiplication of complex numbers.
Fix a surjective homomorphism $s:\IT\to\IQ$.
For every real number $a$ consider the subsets $$\IT^>_{a}:=\{1\}\cup \{x\in\IT:s(x)>a\}\mbox{ \ and \ }\IT_{a}^\ge:=\{1\}\cup \{x\in\IT:s(x)\ge a\}$$and the familes of sets $$\mathcal B'_1=\{\IT^>_a:a<0\},\;\;\mathcal B'_2=\{\IT^\ge_0\},\;\;
\mathcal B'_3=\{\IT^>_0\}\mbox{ \ and \ }\mathcal B'_4=\{\IT^>_a:a>0\}.$$
It is easy to check that  for every $i\in\{1,2,3,4\}$ the family $\mathcal B'_i$ is a base at the identity of a first-countable semigroup
topology on the group $\IT$.  Consequently, for every $i$ the family
$\mathcal B_i=\{S\cap U:S\in \mathcal B'_i,\;1\in U\in\tau\}$ is
a base at the identity of a first-countable functionally Hausdorff semigroup
topology on the group $\IT$. Denote this topology by $\tau_i$ and put $G_i:=(\IT,\tau_i)$.
Clearly, $\tau\subset\tau_1\subset\tau_2\subset\tau_3\subset\tau_4$.
Fix any element $d$ with $s(d)>0$ and put $D=\{nd:n\in\mathbb N\}$.
It is easy to check that the set $D$ is dense in $G_4$.
Let $\{x_n\}_{n\in\omega}$ be any sequence of points of $D$. Since
the group $(G,\tau)$ is sequentially compact, there exists a point $x\in\IT$ and
a subsequence $\{x_{n(k)}\}_{k\in\omega}$ converging to $x$ in $(\IT,\tau)$.
Taking a subsequence, if necessary, we can assume that $s(x_{n(k+1)})\ge s(x_{n(k)})+s(d)$
for each $k$. For any $a>0$ there exists $K>0$ such that $s(x_{n(k)})>s(x)+a$ and so
$x_{n(k)}\in x+\IT^>_a$ for each $k>K$. It follows that $\{x_{n(k)}\}$ converges to $x$ in $G_4$.
It is easy to check that the group $(\mathbb T,\tau)$ is pseudobounded,
the group $G_1$ is $\omega$-pseudo\-bounded, and the groups $G_2$, $G_3$, and $G_4$ are not
$\omega$-pseudobounded.
Observe that for $i\in\{1,2\}$ the preimege $s^{-1}(0)$ is a subgroup of countable index in $\IT$ and
$\tau_i{\restriction}s^{-1}(0)=\tau{\restriction}s^{-1}(0)$.
Since the space $(\IT,\tau)$ is second countable, the paratopological group $G_i$ has a countable
network. Since the space $G_i$ is first-countable, by~\cite[Proposition 2.3]{Rav}, the paratopological group $G_i$ is second countable.
We claim that for $i\in\{2,3\}$ the group $G_i$ does not have a countable network consisting of closed subsets. To derive a contradiction, assume that $\mathcal N$ is a countable family of closed subsets of $G_i$ such that for any $x\in G_i$ and neighborhood $O_x\subseteq G_i$ of $x$ there exists a set $N\in\mathcal N$ with $x\in N\subseteq O_x$.
In particular, for any $x\in \IT$ there exists a closed set $N_x\in\mathcal N$ such that $x\in N_x\subseteq x+\IT^{\ge}_0$. We claim that the set $N_x$ is nowhere dense in $\IT$.
To derive a contradiction, assume that the $\tau$-closure of $N_x$ has nonempty interior in $\IT$.
Since the set $\{nd:n<s(x)\}$ is dense in $\IT$, there is a point $y\in \overline{N_x}^\tau$ such
that $s(y)<s(x)$.  For every neighborhood $V\in\tau$ of $y$ there exists a point $z\in V\cap
N_x\subseteq N_x\subseteq x+\IT^{\ge}_0$ that has $s(z)\ge s(x)>s(y)$. Then $z\in N_x\cap
V\cap(y+\IT^{>}_0)$ and hence $y\in \overline{N_x}^{\tau_i}=N_x\subseteq x+\IT^\ge_0$ and $s(y)\ge s(x)$,
which contradicts the choice of $y$. This contradiction shows that the set $N_x$ is nowhere dense
in $\IT$. Then the union $\IT=\bigcup_{x\in X}N_x$ is a meager subset in $\IT$, which contradicts the
Baireness of the compact space $\IT$. This contradiction proves that the space $G_i$ has no
countable network consisting of closed subsets.
\smallskip
Observe that the set $\IT\bs \IT_0^>$ is closed and discrete in $G_4$. Thus the space $G_4$ has extent $e(G_4)=\mathfrak c$.
The set $\{(x,-x):x\in \IT\}$ is a closed discrete subset of $G_3\times G_3$, which implies that the space $G_3\times G_3$ has extent $e(G_3\times G_3)=\mathfrak c$.

The values of the extent $e(G_3)$ and the Lindel\"of number $l(G_3)$ of the space $G_3$ depend on the homomorphism $s$.
\begin{claim}\label{OstapEx1N} There exists a surjective homomorphism $s:\IT\to\IQ$ such that the set $s^{-1}(0)$
contains a compact subset $K$ of $\IT$ of cardinality $|K|=\mathfrak c$. In this case $e(G_3)=l(G_3)=\mathfrak c$.
\end{claim}

\begin{proof} By~\cite[19.2]{Kech}, the real line contains a Cantor set $C\subseteq\IR$, which is linearly independent over the field $\IQ$. Replacing $C$ by $\frac{2\pi}cC$ where $c\in C$ is any point, we can assume that $2\pi \in C$. Replacing $C$ by a suitable closed subset, we can additionally assume that the $\IQ$-linear hull of $C$ does not contain some real number $r$. Then the set $C\cup\{r\}$ is linearly independent over $\IQ$ and by Zorn's Lemma can be enlarged to a maximal linearly independent set $B$. Let $f:\mathbb R\to \IQ$ be a unique $\IQ$-linear functional such that $f(r)=1$ and $f(B\setminus\{r\})=\{0\}$. Consider the homomorphism $h:\mathbb R\to\IT$, $h:x\mapsto e^{ix}$, and observe that $h^{-1}(1)=2\pi\mathbb Z\subseteq f^{-1}(0)$. Then there exists a unique homomorphism $s:\IT\to\IQ$ such that $f=s\circ h$. The surjectivity of $f$ implies the surjectivity of $s$. Since $h$ has countable preimages of points, the image $K=h(C)$ is an uncountable compact subset of $\IT$ with $s(K)=s(h(C))=f(C)=\{0\}$. By \cite[3.3.2 and 6.5]{Kech}, $|K|=\mathfrak c$.
Since $\tau\subseteq\tau_3$, the set $K$ is a closed in $G_3$.
Since $(x+\IT^>_0)\cap K=\{x\}$ for every $x\in K$, the set $K$ is discrete in $G_3$.
Consequently $\mathfrak c\le e(G_3)\le l(G_3)=\mathfrak c$.
\end{proof}

\begin{claim}\label{OstapEx3N} There is a surjective homomorphism $s:\IT\to\IQ$ such that $\inf s(F)=-\infty$ for every uncountable closed subset $F\subseteq\IT$. In this case $e(G_3)=l(G_3)=\omega$.
\end{claim}
\begin{proof} The construction of $s$ resembles the construction of a Bernstein set.
Since the family $\mathcal F$ of uncountable closed subsets of the real line has cardinality continuum, it can be enumerated as $\mathcal F=\{F_\alpha\}_{\alpha<\mathfrak c}$ so that $F_{\alpha+n}=F_\alpha$ for every limit ordinal $\alpha<\mathfrak c$ and every $n\in\omega$. For every $\alpha<\mathfrak c$ we can  inductively choose a point $x_\alpha\in F_\alpha$ that does not belong to the $\IQ$-linear hull $L_{<\alpha}$ of the set $X_{<\alpha}:=\{2\pi\}\cup\{x_\beta:\beta<\alpha\}$. The choice of the point $x_\alpha$ is always possible since $|F_\alpha|=\mathfrak c>|\w+\alpha|=|L_{<\alpha}|$.
After completing the inductive construction, we obtain a linearly independent subset $\{2\pi\}\cup\{x_\alpha\}_{\alpha<\mathfrak c}$ in the linear space $\IR$ over the field  $\IQ$. Choose a $\IQ$-linear functional
$f:\mathbb R\to\IQ$ such that $f(2\pi)=0$ and $f(x_{\alpha+n})=-n$ for any limit ordinal $\alpha<\mathfrak c$ and any $n\in\omega$. Consider the continuous homomorphism $h:\mathbb R\to \IT$, $h:t\mapsto e^{i t}$, and observe that $h^{-1}(1)=2\pi\cdot \mathbb Z\subseteq f^{-1}(0)$. Consequently, there exists a homomorphism $s:\IT\to\IQ$ such that $f=s\circ h$. We claim that $\inf s(F)=-\infty$ for any uncountable closed set $F\subseteq \IT$. Indeed, the preimage $h^{-1}(F)$ is an uncountable closed subset of $\mathbb R$ and hence $h^{-1}(F)=F_\beta$ for some limit ordinal $\beta$. The choice of the enumeration $\{F_\alpha\}_{\alpha<\mathfrak c}$ ensures that $F_\beta=F_{\beta+n}$ for any $n\in\omega$. By the definition of the functional $f$, we get $s(h(x_{\beta+n}))=f(x_{\beta+n})=-n$. Consequently, $$\inf s(F)=\inf f(h^{-1}(F))=\inf f(F_\beta)\le\inf \{f(x_{\beta+n}):n\in\omega\}=\inf\{-n:n\in\omega\}=-\infty.$$
Now we prove that the space $G_3$ is Lindel\"of. Given any open cover  $\mathcal U$ of the space
$G_3$, for every point $x\in G_3$ choose a neighborhood $O_x\in\tau$ of $x$ and a set
$U_x\in\mathcal U$ such that $(x+\IT_0^>)\cap O_x\subseteq U_x$. For any rational number $a$, consider
the subsets $A_a:=\bigcup\{O_x:x\in\IT\bs\IT^>_a\}$ and $B_a:=\IT\bs A_a\subseteq\IT_a^>$. Since
$B_a\subseteq\IT_a^>$ is a closed subset of $\IT$ with $\inf s(B_a)\ge \inf s(\IT_a^>)=a$, the choice of
the homomorphism $s$ ensures that the set $B_a$ is at most countable.
Since the metrizable compact space $(\IT,\tau)$ is hereditarily Lindel\"of, there exists a countable
subset $C_a$ of $\IT\bs\IT_a^>$ such that $A_a=\bigcup\{O_x:x\in C_a\}$. Consider the countable
subfamily $\mathcal V:=\bigcup_{a\in \IQ}\{U_x:x\in B_a\cup C_a\}$ of $\mathcal U$. We claim that
$\IT=\bigcup\mathcal V$. Indeed, given any point $y\in\IT$, choose a negative rational number
$a<s(y)$. If $y\in B_a$, then $y\in U_y\in \mathcal V$. If $y\notin B_a$, then $y\in
A_a=\bigcup_{x\in C_a}O_x$ and hence $y\in O_x$ for some $x\in C_a\subseteq\IT\setminus\IT_a^>$. It
follows that $s(x)\le a<s(y)$ and hence $y\in O_x\cap(x+\IT_0^>)\subseteq U_x\in\mathcal V$. Therefore
the space $G_3$ is Lindel\"of and has Lindel\"of number $l(G_3)\le\omega$. Since $\omega\le
e(G_3)\le l(G_3)\le \omega$ we have $e(G_3)=l(G_3)=\omega$.
\end{proof}
\end{proof}

In \cite[Theorem~1]{SanMTka} Tkachenko and Sanchis
constructed a Hausdorff $2$-pseudo\-compact Fr\'echet-Urysohn
paratopological group which is not a topological
group, thus answering a second author's
question from an older version of this manuscript.
The next proposition was motivated by a question of Gutik whether the closure of a subgroup of a
countably compact paratopological group is a group.
Also it provides a partial solution to Problem 4.4 from~\cite{TkaTom}: 
Let $G$ be a Hausdorff  feebly  compact paratopological  group  such  that  the  closure  of  every
subgroup of $G$ is a subgroup. Is $G$ a topological group?
\begin{proposition} Let $G$ be a countably compact paratopological group such that the closure
of each cyclic subgroup of $G$ is a group. Then $G$ is a topological group.
\end{proposition}
\begin{proof} Let $x$ be an arbitrary point of $G$ and $H$ be the closure of the cyclic subgroup
$\langle x\rangle$ generated by $x$. By our assumption, $H$ is a subgroup of $G$. We claim that the
group $H$ is left $\omega$-precompact. It suffices to show that $H\subseteq \langle x\rangle\cdot U$
for any neighborhood $U\subseteq H$ of the identity $x^0$. For every point $y\in H$ the density of
the subgroup $\langle x\rangle$ in $H$ yields a point $z\in\langle x\rangle\cap Uy^{-1}$. Then
$y\in z^{-1}U\subseteq \langle x\rangle\cdot U$. The space $H$ is countably compact being a closed subspace
of a countable compact space $G$. Therefore the group $H$ is $2$-pseudocompact. By
Lemma~\ref{2PseBai} $H$ is Baire and by Proposition~\ref{2PseP} $H$ is feebly compact. Thus by Proposition~\ref{PseTG}
$H$ is topologically periodic, and so is the paratopological group $G$. By \cite{BokGur}, $G$ is a topological group.
\end{proof}
The following proposition answers a question of Guran.
\begin{proposition}\label{T2SigCPsTG} Each Hausdorff $\sigma$-compact
feebly compact paratopological group is a compact topological group.
\end{proposition}
\begin{proof} Let $(G,\tau)$ be such a paratopological group and $(G,\tau_{sr})$ be its
semiregularization. By Lemma~\ref{l:semiregularization}, $(G,\tau_{sr})$ is a feebly compact Hausdorff
regular paratopological group and by Proposition~\ref{PseTG}, $(G,\tau_{sr})$ is a pseudocompact
topological group. Write the $\sigma$-compact space $G$ as the countable union
$G=\bigcup_{n\in\omega} K_n$ of compact subsets $K_n$, which remain compact in the topological group
$(G,\tau_{sr})$. Since the pseudocompact topological group $(G,\tau_{sr})$ is $2$-pseudocompact,
it is Baire by Lemma~\ref{2PseBai}, so for some $n\in\omega$ the compact set $K_n$ has nonempty interior
in $(G,\tau_{sr})$ and also in $(G,\tau)$. This means that the space $(G,\tau)$ is locally compact.
By Ellis' Theorem (see Subsection~\ref{subsec:PrelimAC}) $(G,\tau)$ is a topological group.  Being feebly compact, the locally compact
(and hence Ra\u\i kov complete) topological group $(G,\tau)$ is  compact.
\end{proof}

\section{Preservation of feebly compactness by some operations over paratopological groups}\label{CCSProd}
In this section we prove that the feeble compactness of paratopological groups is preserved by Tychonoff products and extensions.
\subsection{Tychonoff products of feebly compact paratopological groups.}
By Tychonoff's theorem, the Tychonoff product of any family of compact spaces is compact.
The productivity of other compact-like properties can be a non-trivial problem.
For instance, the product of a countable family of sequentially
compact spaces is sequentially compact~\cite[Theorem~3.10.35]{Eng}.
But the Cantor cube $D^\mathfrak c$ is not sequentially compact
(see~\cite{Eng}, the paragraph after Example 3.10.38).
On the other hand some compact-like spaces are also preserved by products, see
~\cite[$\S$ 3--4, 7]{Vau} (especially Theorem 3.3, Proposition 3,4, Example 3.15,
Theorem 4.7, and Example 4.15), ~\cite[$\S$ 5]{Ste}, and
~\cite[2.3]{GutRav}.
In~\cite{Tka2},~\cite{Tka3}, and ~\cite{Tka6} Tkachenko considered the productivity of various properties of topological groups.
All but one non-productive properties of spaces listed in~\cite[p.~2]{Tka6},
are (possible, in models) also non-productive for topological groups.
For instance, Nov\'ak~\cite{Nov} and Teresaka~\cite{Ter}  constructed examples of two
countably  compact spaces whose  product is not pseudocompact, see also~\cite[3.10.19]{Eng}.
Using Martin's Axiom, van Douwen~\cite{vDou} constructed
countably compact topological groups $G,H$  whose product $G\times H$ is not countably compact.
Later, under $MA_{countable}$ Hart and van Mill~\cite{HarvMil}
constructed a countably compact
topological group $G$ whose square $G\times G$ is not countably compact.
As far as the authors know, there is no ZFC-example of a family $\mathcal\{G_\alpha:\alpha\in A\}$
of countably compact topological groups whose Tychonoff product $\prod \mathcal\{G_\alpha:\alpha\in A\}$
is not countably compact, see the surveys \cite{Com} and \cite{ComHofRem}.
Feeble compactness of products of spaces is preserved in particular cases, see~\cite{ScaSto}.
Also Dow et al. in Theorem 4.1 of~\cite{DPSW} proved that
the Tychonoff product of a family of sequentially feebly compact 
 spaces is again sequentially feebly compact, and in Theorem 4.3 that every product of feebly compact
spaces, all but one of which are sequentially feebly compact, is feebly compact.
Comfort and Ross in~\cite[Theorem~1.4]{ComRos}
proved that the product of any family of pseudocompact topological
groups is pseudocompact.
The main result of this subsection is a generalization of this result
to paratopological groups.

\begin{theorem}\label{ProdPCompact} The Tychonoff product of an arbitrary family of feebly compact paratopological groups is feebly compact.
\end{theorem}

\begin{proof} Let $(G,\tau)$ be the Tychonoff product of a family
$\{(G_\alpha,\tau_\alpha):\alpha\in A\}$ of feebly compact paratopological groups.
It is easy to see that the semiregularization $(G,\tau_{sr})$ of the space $(G,\tau)$ coincides with the Tychonoff product $\prod_{\alpha\in A}(G,\tau_{\alpha sr})$ of the semiregularizations of the spaces $(G_\alpha,\tau_\alpha)$.
 Let $\tau_{sre}$ be the family of all open neighborhoods of the identity in the paratopological group
 $(G,\tau_{sr})$. It is easy to see that the intersection $H:=\bigcap_{U\in\tau_{sre}}U\cap U^{-1}$ is a normal subgroup of $G$. Let $\pi:G\to G/H$ be the quotient map.
For any element $\alpha\in A$ let $\tau_{\alpha sre}$ be the family of all open neighborhoods of the identity $e_\alpha$ of the  paratopological group
$(G_\alpha,\tau_{\alpha sr})$, $H_\alpha=\bigcap_{U\in\tau_{\alpha sre}}U\cap U^{-1}$,
$\pi_\alpha:(G_\alpha,\tau_{\alpha sr}) \to (G_\alpha,\tau_{\alpha sr})/H_\alpha$ be the quotient map,
and $p_\alpha:(G,\tau_{sr})\to (G_\alpha, \tau_{\alpha sr})$ be the coordinate projection.
Let $\pi:\prod_{\alpha\in A}(G_\alpha,\tau_{\alpha sr})\to \prod_{\alpha\in A} (G_\alpha,\tau_{\alpha sr})/H_\alpha$
be the product of the family of quotient maps $\{\pi_\alpha\}_{\alpha\in A}$. By Proposition 2.3.29 from~\cite{Eng}, the map $\pi$ is open. Since $H=\ol{\{e\}}^{\tau_{sr}}\cap (\ol{\{e\}}^{\tau_{sr}})^{-1}$ and
$H_\alpha=\ol{\{e_\alpha\}}^{\tau_{\alpha sr}}\cap (\ol{\{e\}}^{\tau_{\alpha sr}})^{-1}$ for each $\alpha\in
A$, we see that $H=\prod_{\alpha\in A} H_\alpha=\bigcap_{\alpha\in A} \pi_\alpha^{-1}(H_\alpha)=\ker \pi$.
Let $i:(G,\tau_{sr})/H\to \prod_{\alpha\in A} (G_\alpha,\tau_{\alpha sr})/H_\alpha$ be a map such that $i(xH)=\pi(x)$ for each $x\in G$. Theorem on continuous epimorphisms \cite[p.~42]{Rav} implies that the
map $i$ is well-defined and is a topological isomorphism.
Let $\alpha\in A$ be any element. Then $(G_\alpha,\tau_{\alpha sr})$ is a regular
feebly compact paratopological group (see the remark after Lemma~\ref{SatQua}). So $(G_\alpha,\tau_{\alpha sr})$
is a feebly compact topological group by Lemma~\ref{ArhRezT3}. Then the group $(G_\alpha,\tau_{\alpha
sr})/H_\alpha$ is a pseudocompact topological group. Therefore by Theorem 1.4 from
\cite{ComRos}, the product $\prod_{\alpha\in A} (G_\alpha,\tau_{\alpha sr})/H_\alpha$ is
pseudocompact and so is its isomorphic copy $(G,\tau_{sr})/H$. Consider the quotient map $q:(G,\tau_{sr})\to (G,\tau_{sr})/H$. By Lemma~\ref{l:Tka-T0}, each open subset $U\in\tau_{sr}$ is equal to $q^{-1}(q(U))$. This fact and the feeble compactness of the quotient group $(G,\tau_{sr})/H$ implies the feeble compactness of the topological group $(G,\tau_{sr})$.
Finally, by Lemma~\ref{RegPCompact}, the paratopological group $(G,\tau)$ is
feebly compact.
\end{proof}
\begin{problem} Is the countably pracompactness or $2$-pseudocompactness of paratopological
groups preserved by Tychonoff products?
\end{problem}
This problem is open since 2010, from the first version of this manuscript. Recently it obtained
partial negative answers. Namely, Garc\'ia-Ferreira and Tomita in~\cite{G-FT}
under {\it CH} constructed a countably compact topological group whose square is not selectively
pseudocompact and therefore not countably pracompact. Bardyla, Ravsky, and Zdomskyy~\cite{BRZ}
constructed under {\it MA} a Boolean Tychonoff countably compact topological group without non-trivial convergent
sequences whose square is not countably pracompact.
\subsection{Extensions of feebly compact paratopological groups}\label{CCSExt}
A (topological, algebraic, or topological-algebraic) property $\mathcal P$ of a
semitopological group is called a {\it three-space property} provided
whenever $N$ is a closed normal subgroup of a group $G$ and both $N$ and $G/N$ have
$\mathcal P$, the group $G$ also has $\mathcal P$.
Various three-space properties in the class of topological groups are investigated for
almost a century. Their list is quite long and includes compactness~\cite{Fre},
local compactness~\cite{Ser}, Ra\u\i kov completeness, precompactness, pseudocompactness~\cite{ComRob},
connectedness, and metrizability~\cite{Gra}, see also~\cite{ComRob},~\cite{DieRol}, and~\cite{HewRos}.
On the other hand, in~\cite{BruTka2} Bruguera and Tkachenko constructed an Abelian topological group $G$
and a closed subgroup $N$ of $G$ such that the groups $N$ and $G/N$ are countably compact,
but $G$ is not. Also they studied a lot of three-space properties for compact-like sets in~\cite{BruTka}.
Lin, Lin, and Xie in~\cite{LinLinXie} considered some specific three-space properties.
Recently three-space properties were investigated in some classes of paratopological and semitopological
groups, and not all of three-space properties of topological groups remain three-space also for these wider classes
of groups. Among such properties let us mention the metrizability (even for Tychonoff paratopological groups,
~\cite[Example 3.3]{LiMouWan}),
the first- and second-countability ~\cite{FerSan}.
On the other hand, to be a topological group is a three-space property in the class of paratopological groups
~\cite[Lemma 4]{Rav5}. Since by \cite[Proposition 5.5]{Rav3},
a locally compact paratopological group is a topological group,
it follows that compactness and local compactness are three-space properties in the
class of paratopological groups.
Similarly to the proof of Theorem 3 in \cite{LinLin}, we
can show that both the pseudoboundedness and the $\omega$-pseudoboundedness are three-space properties
in the class of topologized groups.
In ~\cite[p.50-51]{Rav3} the second author showed
that for each infinite cardinal $\lambda$ both $\psi(G)\le\lambda$ and
$d(G)\le\lambda$ are three-space properties in the class of paratopological groups.
For more three-space properties in the class of paratopological groups
see~\cite{LinXie}.
The metrizability, first-countability and second-countability also are not three-space properties
in the class of quasitopological groups, see~\cite{LinLinTan}.
Some three-space properties of semitopological groups are also considered in~\cite{GuoPen} and
~\cite{LiMouWan2}.
The main result of this subsection is Theorem~\ref{PseudoQuot} establishing the three-space property
of the feeble compactness in the class of paratopological groups.
It generalizes the following result
(see Theorem 6.3(c) from the paper~\cite{ComRob} by Comfort and Robertson).
\begin{theorem}\label{t:Tka-3} Let $H$ be a closed normal subgroup of a Hausdorff topological group $G$. If the topological groups $H$ and $G/H$ are feebly compact, then so is the topological group $G$.
\end{theorem}


To generalize Theorem~\ref{t:Tka-3} to the class of paratopological groups, we need to make some preliminary job. First we observe that the closedness of the subgroup $H$ can be easily removed from Theorem~\ref{t:Tka-3}.
\begin{lemma}\label{PseudoQuotTG} Let $H$ be a normal subgroup of a Hausdorff topological group $G$. If the topological groups $H$ and $G/H$ are feebly compact, then so is the topological group $G$.
\end{lemma}
\begin{proof} Observe that the closure $\bar{H}$ of the group $H$ in $G$ is feebly compact (since it contains a dense feebly compact subspace). By the continuity of the homomorphism $h:G/H\to G/\bar H$, $h:gH\mapsto g\bar H$, the quotient group $G/\bar H$ is feebly compact, being a continuous image of the feebly compact space $G/H$. Applying Theorem~\ref{t:Tka-3}, we conclude that the topological group $G$ is feebly compact.
\end{proof}

\begin{lemma}\label{T3QComp_} Let $G$ be a regular paratopological group and $H$ be a compact normal
subgroup in $G$. Then
\begin{enumerate}
\item the quotient paratopological group $G/H$ is regular;
\item if $G/H$ is a topological group, then $G$ is a topological group.
\end{enumerate}
\end{lemma}
\begin{proof} 1. To prove the regularity of $G/H$, it suffices to show that every open neighborhood
$U\subseteq G/H$ of the identity contains the closed neighborhood of the identity. Let $\pi:G\to G/H$ be
the quotient homomorphism. Then $\pi^{-1}(U)=\pi^{-1}(U)H$ is
an open neighborhood of the identity. By the regularity of $G$, it contains a closed
neighborhood $\ol{V}$ of the identity. Then $\ol{V}H\subseteq \pi^{-1}(U)H$.
Propositions 1.12 and 1.13 from~\cite{Rav} imply that the map $\pi$ is open and closed.
Then $\pi(\ol{V}H)=\pi(\ol{V})\subseteq U$ is a closed neighborhood of the identity in $G/H$,
which follows that the paratopological group $G/H$ is regular.
\smallskip
2. By Lemma 5.4 from~\cite{Rav3} or its counterpart in English in~\cite{Rav6},
or by Proposition~\ref{TCCTG}, $H$ is a topological group. Now by
Lemma 4 from~\cite{Rav5}, $G$ is a topological group.
\end{proof}
\begin{corollary}\label{T0Reg} If $G$ is a regular paratopological group
then its $T_0$-reflexion $T_0G$ is regular too.
\end{corollary}

\begin{lemma}\label{l:3reg} Let $(G,\tau)$ be a regular Hausdorff paratopological group  and $H$ be a  normal subgroup of $G$. If the spaces $H$ and $G/H$ are feebly compact, then so is the space $G$.
\end{lemma}
\begin{proof} Let $\breve G=(\breve G,\breve\tau)$ be the Ra\u\i kov completion of the paratopological group $(G,\tau)$ and $\breve G_{sr}=(\breve G,\breve\tau_{sr})$ be the semiregularization of the paratopological group $\breve G$. By Theorem~\ref{t:Raikov}, $\tau=\breve\tau_{sr}{\restriction}
G$ and the group coreflexion of $\breve G_{sr}$ coincides with the Ra\u\i kov completion $\breve G^\sharp$ of the group coreflexion $G^\sharp$ of the paratopological group $G$.
By Proposition~\ref{PseTG}, the regular feebly compact paratopological group $(H,\tau{\restriction}H)$ is a topological group. Then
$\tau^\sharp{\restriction}H=\tau{\restriction}H$ and the topological group $(H,\tau^\sharp{\restriction}H)$ is feebly compact and precompact by Lemma~\ref{PseudoTGPre}. Then its closure $\bar H$ in the Ra\u\i kov-complete topological group $\breve G^\sharp$ is compact, see  \cite[3.7.9]{ArhTka}. Since the identity map $\breve G^\sharp\to\breve G_{sr}$ is continuous and the paratopological group $\breve G_{sr}$ is Hausdorff (by Theorem~\ref{t:Raikov}), the subgroup $\bar{H}$ is compact and closed in $\breve G_{sr}$. The normality of the subgroup $H$ in $G$
implies the normality of the group $\bar{H}$ in the topological group $\breve G^\sharp$. Then the set $\Gamma=G\bar{H}$ is a subgroup of $\breve G_{sr}$.
Consider the quotient paratopological group $\Gamma/\bar H$ (which is regular by Lemma~\ref{T3QComp_}), and  the surjective  homomorphism $\pi:G/H\to \Gamma/\bar H$ assigning to each coset $gH\in G/H$ the coset $g\bar H\in\Gamma/\bar H$.  The continuity of the identity embedding $G\to G\bar H=\Gamma\subseteq \breve G_{sr}$ implies the continuity of the map $\pi:G/H\to \Gamma/\bar H$. Then the space $\Gamma/\bar H$ is feebly compact, being a continuous image of the feebly compact space $G/H$. By Proposition~\ref{PseTG}, the regular feebly compact paratopological group $\Gamma/\bar H$ is a topological group, and by Lemma~\ref{T3QComp_}(2), the paratopological group $\Gamma$ is a topological group and so is its subgroup $G$ (here we use the fact that $\tau=\breve\tau_{sr}{\restriction}G$ established in Theorem~\ref{t:Raikov}). By Lemma~\ref{PseudoQuotTG}, the topological group $G$ is feebly compact.
\end{proof}

\begin{theorem}\label{PseudoQuot} Let $H$ be a normal subgroup of a
paratopological group $G$. If the paratopological groups $H$ and $G/H$ are feebly compact, then $G$ is feebly compact, too.
\end{theorem}
\begin{proof} Consider the semiregularization $G_{sr}$ of the paratopological group $G$ and its $T_0$-reflextion $\Gamma:=T_0G_{sr}$, which is a Hausdorff regular paratopological group (according to Lemma~\ref{l:semiregularization}). Let $\pi:G\to \Gamma$ be the composition of the identity map $G\to G_{sr}$ and the quotient map $G_{sr}\to\Gamma$. The normality of the subgroup $H$ implies the normality of the subgroup $N:=\pi(H)$ in the group $\Gamma$. Being a continuous image of the feebly compact space $H$, the space $N$ is feebly compact. Consider the homomorphism $h:G/H\to\Gamma/N$ assigning to each coset $gH$ the coset $\pi(gH)=\pi(g)N$. The continuity of the  map $\pi:G\to\Gamma$ implies the continuity of the homomorphism $h$. Then the feeble compactness of the quotient group $G/H$ implies the feeble compactness of the quotient group $\Gamma/N$. By Lemma~\ref{l:3reg}, the paratopological group $\Gamma=T_0G_{sr}$ is feebly compact.  Now Lemma~\ref{l:Tka-T0} implies that the space $G_{sr}$ is feebly compact and Lemma~\ref{RegPCompact} implies that so is the space $G$.
\end{proof}
\begin{remark} Theorem~\ref{PseudoQuot} solves positively Problem 5.3
and the quasiregular case of Problem 5.4 of Tkachenko \cite{Tka5}.
As far as the authors know, Problem 5.7 of ~\cite{Tka5}
whether $2$-pseudocompactness is a three-space property in the class of (Hausdorff, $T_3$) paratopological groups remains open.
\end{remark}
\section{Cone topologies}\label{sec:cone-top}
In this section we systematically study the cone topologies on topological groups. Such topologies yield many pathological examples of paratopological groups, see e.g. \cite{Ban}, \cite{BanRav}, \cite{Rav}, \cite{Rav2} or \cite{Tod}.
The cone topologies $\tau_S$ and $\tau_{\vec S}$ on a paratopological group $(G,\tau)$ are determined by a normal subsemigroup $S\subseteq G$ containing the identity $e$. The normality of $S$ means that $g^{-1}Sg=S$ for any element $g\in G$. Let $Z_S=\bigcap_{x\in S}\{z\in S:xz=zx\}$ be the {\em center} of the semigroup $S$.
By definition, the {\em cone topology} $\tau_S$ on $G$ consists of the sets $U\subseteq G$ such that for every $x\in U$ there exists an open neighborhood $V\in\tau$ of $x$ such that $V\cap xS\subseteq U$.
The {\em shifted cone topology} $\tau_{\vec S}$ on $G$ consists of the sets $U\subseteq G$
such that for every $x\in U$ there exists an open neighborhood $V\in\tau$ of $x$ and a
point $z\in Z_S$ such that $V\cap xzS\subseteq U$.
It can be shown that $(G,\tau_S)$ and $(G,\tau_{\vec S})$ are paratopological groups.
A particular case of this construction appears when $G$ is a real linear space and
$S$ is a cone in $G$, that is a subset of $G$ such that $\lambda x+\mu y\in S$
for any points $x,y\in S$ and non-negative real numbers $\lambda,\mu$,
which explains the term ``cone topology".
The Sorgenfrey topology on the real line $\IR$ is a cone topology $\tau_S$ determined by the cone $S=[0,+\infty)$.
If the topology $\tau$ on $G$ is antidiscrete, then the paratopological groups $(G,\tau_S)$ and $(G,\tau_{\vec S})$ will be denoted by $G_S$ and $G_{\vec S}$, respectively.
Observe that for any semigroup topology $\tau$ on $G$, the identity maps $(G,\tau_S)\to G_S$ and $(G,\tau_{\vec S})\to G_{\vec S}$ are continuous. So if a topological
property $\mathcal P$ is preserved by continuous isomorphisms (for instance, the  compactness, countably compactness,
feeble compactness) and we wish the group $(G,\tau_S)$ or $(G,\tau_{\vec S})$ to have the property $\mathcal P$, then
we must assure the group $G_S$ or $G_{\vec S}$ to have the property $\mathcal P$, too.
So we start with studying the paratopological groups $G_S$ and $G_{\vec S}$.
As expected, in this basic case the cone topologies have very weak separation properties, and
even $T_1$ or $T_2$ cases are degenerated, see Subsections~\ref{Antidicrete1} and \ref{Anti2}. But
non-Hausdorffness of these constructions will be compensated
by a Hausdorff topology $\tau$ in Subsection~\ref{Cowide}.
It turned out that there is a simple interplay between the algebraic properties
of the semigroup $S$ and compact-like properties of the paratopological groups $G_S$ and $G_{\vec S}$, see Subsections~\ref{Antidicrete1} and \ref{Anti2}.
In these subsections  we assume that $G$ is an Abelian group and $S\ni 0$ is a subsemigroup of $G$.
We restrict our attention here to Abelian groups, because
authors' experience suggests that when we check an implication of
topological properties for a paratopological groups, then
positive results are proved for all groups, but negative are illustrated by Abelian
counterexamples.
In the sequel we shall need an elementary lemma.
\begin{lemma}\label{MajorOfWhole}If there is an element $a\in S$ such that $S\subseteq a-S$, then $S$
is a group.
\end{lemma}
\begin{proof} Since $2a\in S\subseteq a-S$, there exists an element $a'\in S$ such that $2a=a-a'$. Thus
$-a=a'\in S$. Since $S$ is a subsemigroup of a group,
it suffices to remark that for each $x\in S$ we have $-x\in -S\subseteq S-a=S+a'\subseteq S$.
\end{proof}
\subsection{The paratopological groups $G_S$}\label{Antidicrete1}
In this section we consider the interplay between algebraic properties
of the subsemigroup $S\ni 0$ of an Abelian group $G$ and topological properties of the paratopological group $G_S$ whose topology consists of the sets $A+S$ where $A\subseteq G$ is an arbitrary subset.
It is easy to see that the closure $\ol S$ of the semigroup $S$ is equal to the clopen subgroup $S-S$ of $G_S$. The following propositions and lemmas can be easily derived from the definition of the topology of the paratopological group $G_S$.
\begin{proposition}\label{G_ST0} The group $G_S$ is $T_0$ iff $S\cap (-S)=\{0\}$.\qed
\end{proposition}
\begin{proposition} The group $G_S$ is $T_1$ iff $S=\{0\}$.\qed
\end{proposition}
\begin{lemma} A subset $K$ of $G_S$ is compact iff $K\subseteq F+S$ for some finite subset $F\subseteq K$.\qed
\end{lemma}
\begin{lemma} For each subset $A\subseteq G$ its closure in $G_S$ coincides with the set $A-S$.\qed
\end{lemma}
It is easy to prove the following characterization.
\begin{proposition}\label{p:top-period}
The paratopological group $G_S$ is topologically periodic iff $S$ is a subgroup of the group $G$ and
the quotient group $G/S$ is periodic.\qed
\end{proposition}
Next, we characterize feebly compact and countably pracompact spaces among paratopological groups $G_S$.
\begin{proposition}\label{G_SPCompact} The following conditions are equivalent:
\begin{enumerate}
\item The space $G_S$ is feebly compact.
\item The space $G_S$ is countably pracompact.
\item The space $G_S$ is sequentially pracompact.
\item The subgroup $S-S$ has finite index in $G$.
\end{enumerate}
\end{proposition}
\begin{proof} The implications $(3)\RT(2)\RT(1)$ are trivial.
\smallskip
$(4)\RT(3)$ Since the subgroup $S-S$ has a finite index in $G$ then there exists a finite set
$F\subseteq G$ such that $F+S-S=G$. Thus $F+S$ is dense in $G_S$. Any sequence of points of
$F+S$ has a subsequence contained in a set $f+S$ for some $f\in F$, so it converges to
$f$.
\smallskip
$(1)\RT (4)$ If the subgroup $H=S-S$ has infinite index in $G$, then we can find an infinite subset
$I\subseteq G$ such that the indexed family $\{xH\}_{x\in I}$ is disjoint.
Since $\{xH\}_{x\in I}$ is an infinite locally finite family of open sets in $G_S$, the space $G_S$
is not feebly compact.
\end{proof}

\begin{proposition}\label{G_SCompact} The following conditions are equivalent:
\begin{enumerate}
\item The space $G_S$ is compact;
\item The space $G_S$ is $\omega$-bounded.
\item The space $G_S$ is totally countably compact.
\item The paratopological group $G_S$ is precompact.
\item The semigroup $S$ is a subgroup of $G$ of finite index.
\end{enumerate}
\end{proposition}
\begin{proof}
The implications $(5)\RT(1)\RT(2)\RT(3)$ are trivial.
\smallskip
$(3)\RT(4)$ If the space $G_S$ is totally countably compact then the set $\ol{\{0\}}=-S$ is compact
and hence is contained in $F+S$ for some finite subset $F$ of $-S$.
Since $G_S$ is feebly compact, by Proposition~\ref{G_SPCompact},
the subgroup $S-S$ has finite index in $G$ and hence $G=D+S-S$ for some finite set $D\subseteq G$.
Thus $G\subseteq D+S-S\subseteq D+S+F+S\subseteq D+F+S$,
which follows that the paratopological group $G_S$ is precompact.
\smallskip

$(4)\RT (5)$ Assuming that the paratopological group $G_S$ is precompact, choose a finite subset $F$ of $G$ such that $F+S=G$. Then  $F+(S-S)=G$ and the subgroup $H=S-S$ has finite index in $G$. For each element $x$ of the finite set $F\cap H\subseteq H=S-S$ find elements $y_x,z_z\in S$ such that $x=y_x-z_x$. Consider the sum $z:=\sum_{x\in F\cap H}z_x\in S$ and observe that $y_x-z_x\in S-z+S=-z+S$ and hence $S-S=(F\cap H)+S\subseteq -z+S+S=-z+S$, which implies that  $S-S=z+S-S\subseteq S$ and $S=S-S$ is a subgroup (of finite index in $G$).
\end{proof}

\begin{proposition}\label{G_S2PCompact} The following conditions are equivalent:
\begin{enumerate}
\item The space $G_S$ is sequentially compact.
\item The space $G_S$ is countably compact.
\item A power $G_S^\kappa$ of the space $G_S$ is countably compact for each cardinal $\kappa$.
\item The space $G_S$ is countably o-compact.
\item The paratopological group $G_S$ is $2$-pseudocompact.
\smallskip
\item The subgroup $S-S$ has finite index in $G$ and each countable subset $C\subseteq S-S$ is contained in $a-S$ for some $a\in S$;
\item The subgroup $S-S$ has finite index in $G$ and each countable subset $C\subseteq S$ is contained in $a-S$ for some $a\in S$;
\item The subgroup $S-S$ has finite index in $G$ and each  infinite subset $C\subseteq S$ has infinite intersection with the set $a-S$ for some $a\in S$;
\item Each  infinite subset $C\subseteq G$ has infinite intersection with the set $a-S$ for some $a\in G$.
\end{enumerate}
\end{proposition}
\begin{proof} We shall prove the implications $(1)\RT(2)\RT(5)\RT(9)\RT(8)\RT(7)\RT(6)\RT(1)$, $(9)\RT(4)\RT(5)$, and $(6)\RT(3)\RT(2)$. Among those implications, $(1)\RT (2)\RT (5)$, $(4)\RT(5)$ and $(3)\RT(2)$ are trivial.
\smallskip
$(5)\RT(9)$ Let $C$ be a countable infinite subset of $G$ and $C=\{x_n\}_{n\in\omega}$ be an enumeration of $C$ such that $x_n\ne x_m$ for any numbers $n<m$. For every
$n\in\omega$ consider the open set $U_n=\bigcup_{i\ge n}(x_i+S)$ of $G_S$. Then $\{U_n:n\in\omega\}$ is a decreasing sequence of nonempty
open subsets of $G_S$ and by the $2$-pseudocompactness of the paratopological group $G_S$, there is a point $b\in G$ such
that $(b+S)\cap -U_n\not=\0$ for each $n\in\omega$. Then for each $n\in\omega$ there is a
number $i\ge n$ such that $(b+S)\cap -(x_i-S)\not=\0$ and hence $-x_i\in b+S$. Let $B$
be the set of all such elements $x_i$. Then $B\subseteq C$ is an infinite subset of $-b-S$.
\smallskip
$(9)\RT(8)$ Assuming that the subgroup $H=S-S$ has infinite index, we can find an infinite subset $I\subseteq G$ that has one-point intersection with each coset $x+H$. Then for every $a\in G$ the intersection $I\cap(a-S)\subseteq I\cap (a+H)$ contains at most one point, which contradicts the condition (9). This contradiction shows that the subgroup $S-S$ has finite index in $G$. For every infinite set $C\subseteq S$ the condition (9) yields a point $a\in G$ such that $C\cap(a-S)$ is infinite and hence $a\in C+S\subseteq S+S$, witnessing that the condition (7) holds.
\smallskip
$(8)\RT(7)$ Let $C=\{c_n\}_{n\in\omega}\subseteq S$ be a countable subset of $S$.
Consider the set $D=\{d_n:n\in\omega\}$, where $d_n=\sum_{i=0}^n c_i$ for each
$n\in\omega$.
If the set $D$ if finite, then $C\subseteq (\sum_{d\in D} d)-S$.
Now suppose the set $D$ is infinite. By (8), there are
an element $a\in S$ and an infinite subset $B$ of $D$ such that $B\subseteq a-S$. Therefore for
each number $i\in\omega$ there exists a number $n\in\omega$ such that $i\le n$ and $d_n\in B$. Then
$c_i\in d_n-S\in a-S-S\in a-S$.
\smallskip
$(7)\RT (6)$. Let $\{c_n:n\in\omega\}$ be a countable subset of the subgroup $S-S$. Fix sequences $\{a_n:n\in\omega\}$ and $\{b_n:n\in\omega\}$ of $S$
such that $c_n=b_n-a_n$ for each $n$. By (7), there exist an element $a\in S$ such that $b_n\in a-S$ for each $n$.
Then $c_n\in a-S$ for each $n$.
\smallskip
$(6)\RT(1)$ To show that the space $G_S$ is sequentially compact, fix any sequence $\{x_n\}_{n\in\omega}$ in $G$.
By (6), the subgroup $S-S$ has finite index in $G$ and hence $G=(S-S)+F$ for some finite subset $F$. By the Pigeonhole Principle, for some $z\in F$ the set $\Omega:=\{n\in\omega:-x_n\in S-S+z\}$ is infinite. By the condition (6),  the set $\{-x_n-z\}_{n\in\Omega}\subseteq S-S$ is contained in the set $a-S$ for some $a\in S$. Then for every $n\in\Omega$ we have $-x_n-z\in a-S$ and hence $x_n\in -a-z+S$, which means that $-a-z$ is the limit of the (convergent) subsequence $(x_n)_{n\in\Omega}$ of the sequence $(x_n)_{n\in\omega}$, and the space $G_S$ is sequentially compact.
\smallskip
$(9)\RT(4)$ To show that the space $G_S$ is countably o-compact, fix a non-increasing sequence
$(U_n)_{n\in\omega}$ of nonempty open sets in $G_S$. For every $n\in\omega$ fix a point $u_n\in
U_n$ and observe that $u_n+S\subseteq U_n$. By the condition (9), there exists a point $a\in G$ such
that the set $\Omega=\{n\in\omega:u_n\in a-S\}$ is infinite. Then $a\in u_n+S\subseteq U_n$ for every
$n\in\Omega$ and hence $a\in \bigcap_{n\in\Omega}U_n=\bigcap_{n\in\omega}U_n$.
\smallskip
$(6)\RT (3)$ Let $\kappa$ be an arbitrary cardinal. By (6) the subgroup $H=S-S$ has finite index in $G$, which implies that the space $G_S$ is homeomorphic to the product $H_S\times D$ of the paratopological group $H_S$ with a finite discrete space $D$.
Then the power $G_G^\kappa$ is homeomorphic to $H_S^\kappa\times D^\kappa$.
By (6), for any sequence $\{x_n\}_{n\in\omega}$ in $H^\kappa$ there exists a point $a\in H^\kappa$ such that $\{x_n\}_{n\in\omega}\subseteq a+S^\kappa$, which implies that $a$ is a limit of the sequence $(x_n)_{n\in\omega}$. Therefore, the space $H_S^\kappa$ is sequentially compact and the product $H_S^\kappa\times D^\kappa$ is countably compact and so is its topological copy $G_S^\kappa$.
\end{proof}

\begin{proposition}\label{G_Sseq-compact} The following statements are equivalent:
\begin{enumerate}
\item Every sequence in $G_S$ is convergent.
\item A power $G_S^\kappa$ of the space $G_S$ is sequentially compact for each cardinal $\kappa$.
\item The power $G_S^{\mathfrak c}$ is sequentially compact.
\smallskip
\item Each countable subset $C\subseteq G$ is contained in $a-S$ for some $a\in S$.
\item $G=S-S$ and each countable subset $C\subseteq S$ is contained in $a-S$ for some $a\in S$.
\item $G=S-S$ and each  infinite subset $C\subseteq S$ has infinite intersection with the set $a-S$ for some $a\in S$.
\item $G=S-S$ and each  infinite subset $C\subseteq G$ has infinite intersection with the set $a-S$ for some $a\in G$.
\end{enumerate}
\end{proposition}
\begin{proof} The implication $(4)\RT(1)\RT(2)\RT(3)$ are trivial and the equivalences $(4)\Leftrightarrow(5)\Leftrightarrow(6)\Leftrightarrow (7)$ follow from the corresponding equivalences in Proposition~\ref{G_S2PCompact}. It remains to prove that $(3)\RT(7)$. Assuming that the power of $G_S^{\mathfrak c}$ is sequentially compact, we first prove that the subgroup $H=S-S$ coincides with the group $G$.   Assuming that $G\ne H$, we conclude that the space $G_S$ is homeomorphic to the product $H_S\times D$ for some discrete space $D$ containing more than one point. The power $D^{\mathfrak c}$ contains a closed topological copy of the Stone-\v Cech compactification $\beta\omega$ and hence is not sequentially compact. Since the space $G_S^{\mathfrak c}$ containst a closed topological copy of the space $D^{\mathfrak c}$, the space $G_S^{\mathfrak c}$ is not sequentially compact, which contradicts (3). This contradiction shows that $G=S-S$. The sequential compactness of $G_S^{\mathfrak c}$ implies the sequential compactness of the space $G_S$. Now Proposition~\ref{G_S2PCompact} ensures that the condition (7) holds.
\end{proof}
Finally, let us observe that the following trivial characterization of the final compactness of the spaces $G_S$.
\begin{proposition}\label{G_SFCompact} The space $G_S$ is Lindel\"of iff the paratopological group $G_S$ is left
$\omega$-precompact.\qed
\end{proposition}
\subsection{The paratopological groups $G_{\vec S}$.}\label{Anti2}
In this section we survey some compact-like properties of the paratopological groups
$G_{\vec S}$. Here $G$ is an Abelian group, $S\ni 0$ is a subsemigroup of $G$ and $G_{\vec S}$ is the group $G$ endowed with the topology $\tau_{\vec S}$ consisting of all sets $U\subseteq G$ such that for every $u\in U$ there exists $x\in S$ such that $u+x+S\subseteq U$. The topology $\tau_{\vec S}$ will be called the
{\em a shift cone topology} generated by the subsemigroup $S$ of $G$.
It is easy to check that the closure $\ol S$ of the semigroup $S$ is equal to the subgroup $S-S$, which is closed-and-open in $G_{\vec G}$. Let us also observe that $\tau_S\subseteq \tau_{\vec S}$.
\begin{proposition} The space $G_{\vec S}$ is Hausdorff iff $S=\{0\}$.\qed
\end{proposition}
Since $\tau_S\subseteq\tau_{\vec S}$, Proposition~\ref{p:top-period} implies the following characterization.
\begin{proposition}
The group $G_{\vec S}$ is topologically periodic iff $S$ is a subgroup of the group $G$ and
the quotient group $G/S$ is periodic.\qed
\end{proposition}

\begin{proposition}\label{G_S*T} The following conditions are equivalent:
\begin{enumerate}
\item The space $G_{\vec S}$ is $T_1$.
\item The space $G_{\vec S}$ is $T_0$.
\item $S=\{0\}$ or $S$ is not a group.
\item $S=\{0\}$ or $\bigcap_{x\in S}(x+S)=\0$.
\end{enumerate}
\end{proposition}
\begin{proof} Since the implications $(1)\RT(2)\RT(3)$ are trivial, it remains to show that $(3)\RT(4)\RT(1)$.
\smallskip

$(3)\RT(4)$ Assuming that (4) does not hold, we conclude that $S\ne\{0\}$ and the intersection $\bigcap_{x\in S}(x+S)$ contains some point $y$. Then $y-S\subseteq S$, $y\in S+S=S$ and $S-S=y-S+S\subseteq S+S=S$, which means that $S=S-S\ne\{0\}$ is a subgroup of $G$. But this contradicts (3).
\smallskip

$(4)\RT(1)$ Assume that the condition (4) holds. To show that $G_{\vec S}$ is a $T_1$-space, fix any non-zero element $x\in G$. If $S=\{0\}$, then $S\in\tau_{\vec G}$ is a neighborhood of $0$ that does not contain $x$. If $\bigcap_{y\in S}(y+S)=\0$, then $x\notin y+S$ for some $y\in S$ and $\{0\}\cup (y+S)\in\tau_{\vec S}$ is a neighborhood of $0$ that does not contain $x$.
\end{proof}

\begin{proposition}\label{G_S*PCompact} The space $G_{\vec S}$ is feebly compact if and only if the subgroup $S-S$ has finite index in $G$.
\end{proposition}
\begin{proof} If the space $G_{\vec S}$ is feebly compact, then so is its continuous image $G_S$. Now Proposition~\ref{G_SPCompact} implies that the subgroup $S-S$ has finite index in $G$.
\smallskip
Now assume that the subgroup $H=S-S$ has finite index in $G$ and find a finite set $F\subseteq G$
such that $G=F+H$. Let $\mathcal U$ be a locally finite family of nonempty open sets in $G$.
Therefore for each $x\in F$ we can pick $s_x\in S$ such that a family $\mathcal U_x=\{U\in\mathcal
U: U\cap (x+s_x+S)\ne\0\}$ is finite. Let $U\in\mathcal U$ be any set. Since $U\subseteq F+S-S$,
there exist $x\in F$ and $s,s'\in H$ such that $x+s-s'\in U$. Since the set $U$ is open, there
exists $s''\in S$ such that $x+s-s'+s''+S\subseteq U$. Then $x+s_x+s+s''\in (x+s_x+S)\cap
(x+s-s'+s''+S)$, so $U\in\mathcal U_x$. Thus the family $\mathcal U=\bigcup_{x\in F}\mathcal U_x$
is finite, which follows that the space $G_{\vec S}$ is feebly compact.
\end{proof}

\begin{proposition}\label{G_S*Compact} The following conditions are equivalent:
\begin{enumerate}
\item The space $G_{\vec S}$ is compact.
\item The space $G_{\vec S}$ is $\omega$-bounded.
\item The space $G_{\vec S}$ is totally countably compact.
\item The space $G_{\vec S}$ is sequentially compact.
\item The space $G_{\vec S}$ is countably compact.
\item The paratopological group $G_{\vec S}$ is precompact.
\item The semigroup $S$ is a subgroup of finite index in $G$.
\end{enumerate}
\end{proposition}
\begin{proof} The implications $(1)\RT (2)\RT (3)$ are trivial.
\smallskip
$(3)\RT(7)$ Assuming that the space $G_{\vec S}$ is totally countably compact, we conclude that its continuous image $G_S$ is totally countably compact, too.
By Proposition~\ref{G_SCompact}, $S$ is a subgroup of finite index in $G$.
The implications $(7)\RT(1)\RT(6)$ are trivial and $(6)\RT(7)$ follows from
 Proposition~\ref{G_SCompact}.
 \smallskip
 The implications $(7)\RT(4)\RT(5)$ are trivial. Finally we prove that $(5)\RT(7)$. Assume that the space $G_{\vec G}$ is countably compact and hence feebly compact. By Proposition~\ref{G_S*PCompact}, the subgroup $S-S$ has finite index in $G$. Next, we check that $S-S=S$. Given any $x\in S$, we should prove that $-x\in S$. This is clear if $x$ has finite order in $G$. So, assume that the order of $x$ is infinite and hence the set $X=\{-nx:n\in\mathbb N\}$ is infinite. By the countable compactness, this set has an accumulation point $a\in G$. Then the neighborhood $a+S\in\tau_{\vec S}$ of $a$ contains some point $-nx$ and hence $-nx\in a+S$ and $a\in -nx-S\subseteq -S$. Since $a\in -S$ is an accumulation point of the set $X$, the neighborhood $\{a\}\cup S=\{a\}\cup (a+(-a)+S)\in\tau_{\vec S}$ of $a$ contains some point $-mx\ne a$ with $m\in\mathbb N$. Then $-x=-mx+(m-1)x\in S+S\in S$ and hence $S$ is a subgroup of $G$.
\end{proof}
\begin{proposition}\label{G_S*CPCompact} The following conditions are equivalent:
\begin{enumerate}
\item The space $G_{\vec S}$ is countably pracompact.
\item The space $G_{\vec S}$ is sequentially pracompact.
\item The subgroup $S-S$ has finite index in $G$ and $S\subseteq C-S$ for some countable set $C\subseteq S$.
\end{enumerate}
\end{proposition}
\begin{proof} The implication $(2)\RT(1)$ is trivial. To prove that $(1)\RT(3)$, assume that the space $G_{\vec S}$ is countably pracompact and hence feebly compact. By Proposition~\ref{G_S*PCompact}, the subgroup $H=S-S$ has finite index in $G$.
 If $S$ is countable, then we can put $C=S$ and conclude that $S\subseteq C-S$. So, assume that $S$ is uncountable.
 The countable pracompactness of the paratopological group $G_{\vec S}$ implies the countable pracompactness of its clopen subgroup $H_{\vec S}$. Consequently, $H_{\vec S}$ is countably compact at some dense subset $A\subseteq S-S$.
 If $A$ is finite, then choose any point $b\in S\setminus A$ and observe that for every $s\in S$ the basic neighborhood $\{b\}\cup(b+s+S)$ of $b$ intersects $A$ and hence $s\in -b+A-S\subseteq A-S$. Since $A\subseteq S-S$ is finite, there exists a finite set $C\subseteq S$ such that $A\subseteq C-S$. Then $S\subseteq A-S\subseteq C-S-S\subseteq C-S$ and we are done.
Now consider the case of infinite set $A$. In this case we can choose a countable infinite subset $B\subseteq A$ and using the countable compactness of $H_{\vec S}$ at $A$, find an accumulation point $b\in H$ of the set $B$. For every $s\in S$ the basic neighborhood $\{b\}\cup(b+s+S)$ of $b$ has infinite intersection with the set $B$. Write $b$ as $b=b_+-b_-$ for some $b_+,b_-\in S$ and observe that the set $s+S\supset b_++s+S=b+s+b_-+S$ has infinite intersection with $B$. Consequently, $s\in B-S$ for every $s\in S$.
Since $B\subseteq A\subseteq S-S$, we can find a countable set $C\subseteq S$ such that $B\subseteq C-S$. Then $S\subseteq B-S\subseteq C-S-S\subseteq C-S$ and we are done.
\smallskip
$(3)\RT(2)$ Assume that the condition (3) holds. Then the subgroup $H=S-S$ has finite index in $G$
and the space $G_{\vec S}$ is homeomorphic to the product $H_{\vec S}\times F$ for some finite discrete
space $F$. The condition (3) yields a countable set $C=\{c_n\}_{n\in\omega}\subseteq S$ such that
$S\subseteq C-S$ and hence $H=S-S\subseteq C-S-S=C-S$. Consider the countable set
$K:=\{0\}\cup \{\sum_{i=0}^n c_i:n\in\omega\}$. We claim that $K$ is a sequence converging to zero
and a dense set in $H_{\vec S}$. To see the latter it suffices to check that for any elements
$x,y,z\in S$ the set $x-y+z+S$ intersects $K$. It follows from $S-S\subseteq C-S$ that $x-y+z\in c_n-S$
for some $n\in\omega$. Then
$$K\ni\sum_{i=0}^nc_i=c_n+\sum_{i=0}^{n-1}c_i\in (x-y+z+S)+\sum_{i=0}^{n-1}c_i
\subseteq (x-y+z)+S+S=x-y+z+S,$$ witnessing that $K$ is dense in $H_{\vec S}$.
To see that the sequence $\{\sum_{i=0}^nc_i\}_{n\in\omega}$ converges to zero, take any $s\in S$ and consider the basic
neighborhood  $\{0\}\cup(s+S)\in\tau_{\vec S}$ of zero. It follows from $s\in S\subseteq C-S$ that
$s\in c_k-S$ for some $k\in\omega$. Then for every $m\ge k$ we get
$\sum_{i=0}^mc_i\in c_k+S\subseteq s+S$, which means that the neighborhood $\{0\}\cup(s+S)$ contains
all but finitely many points of the set $K$. 
\end{proof}

\begin{proposition}\label{G_S*2PCompact} The following conditions are equivalent:
\begin{enumerate}
\item The paratopological group $G_{\vec S}$ is $2$-pseudocompact.
\item Each power of $G_{\vec S}$ is $2$-pseudocompact.
\item The space $G_{\vec S}$ is countably o-compact.
\item The subgroup $S-S$ has finite index in $G$ and each countable subset of $S-S$ is contained in the set $a-S$ for some $a\in S$.
\item Each infinite subset of $G$ has infinite intersection with the set $a-S$ for some $a\in G$.
\end{enumerate}
\end{proposition}
\begin{proof} The implications $(3)\RT(1)$ is trivial and $(1)\RT(4)\RT(5)$ follows from  Proposition~\ref{G_S2PCompact} (and the preservation of $2$-pseudocompactness by continuous homomorphisms of paratopological groups).
\smallskip
$(5)\RT(3)$ Assume that the condition (5) is satisfied. To prove that the space $G_{\vec S}$ is
countably o-compact, fix a decreasing sequence $(U_n)_{n\in\omega}$ of nonempty open sets in
$G_{\vec S}$. For every $n\in \omega$ fix a point $u_n\in U_n$
such that $u_n+S\subseteq U_n$. By the condition (5), there exists a point $a\in G$
such that the set $\Omega=\{n\in\omega:u_n\in a-S\}$ is infinite.
Then for every $n\in\Omega$ the point $a$ belongs to the set $u_n+S\subseteq U_n$ and
hence $a\in\bigcap_{n\in\Omega}U_n=\bigcap_{n\in\omega}U_n$,
witnessing that the space $G_{\vec S}$ is countably o-compact.
\smallskip
Since $(2)\RT(1)$, it remains to prove that $(4)\RT(2)$. Given any cardinal $\kappa$, fix a
decreasing sequence $(U_n)_{n\in\omega}$ of nonempty open sets in the space $G^\kappa_{\vec S}$.
For every $n\in\omega$ fix an element $u_n\in U_n$ and find a finite set $F_n\subseteq\kappa$ and a
function $s_n\in S^\kappa$ such that the basic open set $$V_n:=\bigcap_{\alpha\in F_n}\{u\in
G_{\vec S}^\kappa:u(\alpha)\in \{u_n(\alpha)\}\cup(u_n(\alpha)+s_n(\alpha)+S)\}$$is contained in
$U_n$. Since the subgroup $H=S-S$ has finite index in $G$, there exists a finite set $D\subseteq G$
such that $H+D=G$ and $D\cap(x+H)=\{x\}$ for every $x\in D$. For every $n\in\omega$ the element
$u_n\in G^\kappa$  can be uniquely written as $u_n=h_n+f_n$ for some $h_n\in H^\kappa$ and $f_n\in
D^\kappa$. By the condition (4), there exists a function $a\in H^\kappa$ such that
$\{h_n(\alpha)+s_n(\alpha):n\in\omega\}\subseteq a(\alpha)-S$ for every $\alpha\in\kappa$. Then for
every $\alpha\in\kappa$ we have $$a(\alpha)+f_n(\alpha)\in
h_n(\alpha)+f_n(\alpha)+s_n(\alpha)+S=u_n(\alpha)+s_n(\alpha)+S$$and hence $a+f_n\in V_n\subset
U_n$ for every $n\in\omega$. Since the space $D^\kappa$ is compact, the sequence
$(f_n)_{n\in\omega}$ has an accumulation point $f\in D^\kappa$. We claim that
$-a-f\in\bigcap_{n\in\omega}\overline{-U_n}$. This will follow from the non-decreasing property of
$(U_n)_{n\in\omega}$ as soon as we show that any open neighborhood $W\subseteq G_{\vec S}$ of $-a-f$
intersects infinitely many sets $-U_n$. For the neighborhood $W$ we can find a finite set
$E\subseteq\kappa$ and a function $s\in S^\kappa$ such that $\bigcap_{\alpha\in E}\{w\in G^\kappa:
w(\alpha)=-a(\alpha)-f(\alpha)\}\subseteq W$. Since $f$ is an accumulation point of the sequence
$(f_n)_{n\in\omega}$ in $D^\kappa$, the set $\Omega:=\bigcap_{\alpha\in
E}\{n\in\omega:f_n(\alpha)=f(\alpha)\}$ is infinite. For every $n\in\Omega$  the point $-a-f_n$
belongs to $W$. On the other hand, $a+f_n\in V_n\subseteq U_n$ and hence $-a-f_n\in -U_n$.
\end{proof}
Propositions~\ref{G_S*Compact} and \ref{G_S*2PCompact} imply the following characterization.
\begin{corollary} If the semigroup $S$ is countable, then the space $G_{\vec S}$ is compact if and only if it is $2$-pseudocompact.
\end{corollary}
\begin{proof}  The ``if'' part is trivial. Suppose that the group  $G_{\vec S}$ is $2$-pseudocompact.
It suffices to show that $S$ is a group. By Proposition~\ref{G_S*2PCompact},
$S-S\subseteq a-S$ for some $a\in S$. In particular, $2a=a-a'$ for some $a'\in S$.
Then $a=-a'\in -S$ and so $S-S\subseteq a-S\subseteq -S-S\subseteq -S$.
\end{proof}

\begin{corollary}\label{G_S*CompactTop} The space $G_{\vec S}$ is compact if and only if the paratopological group $G_{\vec S}$ is countably pracompact and
$2$-pseudocompact.
\end{corollary}
\begin{proof} The ``only if'' part is trivial. To prove the ``if'' part, assume that the
paratopological group $G_{\vec S}$ is countably pracompact and
$2$-pseudocompact. By Propositions~\ref{G_S*CPCompact} and\ref{G_S*2PCompact}, the subgroup $H=S-S$
has finite index in $G$, $S\subseteq C- S$ for some countable subset $C\subseteq S$ such that
$C\subseteq a-S$ for some $a\in S$. Then $S\subseteq C-S\subseteq a-S-S=a-S$ and
$-S\subseteq S-S=-a+(S-S)=(-a+S)-S\subseteq -S-S=-S$, which implies that $S=S-S$ is a subgroup of $G$.
By  Proposition~\ref{G_S*Compact}, the space $G_{\vec S}$ is compact.
 \end{proof}

\begin{proposition}\label{G_S*Baire} If space $G_{\vec S}$ is Baire then either $S$ is a group or $S\not\subseteq C-S$ for any countable set $C\subseteq S$.
\end{proposition}
\begin{proof} Suppose that $S\subseteq C-S$ for some countable subset $C\subseteq S$. Then $S-S\subseteq C-S-S=C-S$. Observe that for every $c\in C$ the set $c-S$ is closed in $G_{\vec S}$.
Since the subgroup $S-S$ is open in the Baire space $G_{\vec S}$, there exists $c\in C$ such that
the closed set $c-S$ has nonempty interior in $S-S$. Then $a+S\subseteq c-S$ for some $a\in S-S$ and
hence $(a-c)+S\subseteq -S$. Then $-S\subseteq S-S=(a-c)+S-S\subseteq -S-S=-S$ and hence $S=S-S$ is a
subgroup of $G$.
\end{proof}
\begin{corollary}\label{G_S*BaireCP} The space $G_{\vec S}$ is compact if and only if $G_{\vec S}$ is countably pracompact and Baire.
\end{corollary}
\begin{proof} If the space $G_{\vec S}$ is compact, then by Proposition~\ref{G_S*Compact} , the semigroup $S$ is a subgroup of finite index in $G$. The definition of the topology $\tau_{\vec S}$ implies that each open dense subset of $G$ coincides with $G$, which implies that the space $G_{\vec S}$ is Baire. Being compact, the space $G_{\vec S}$ is countably pracompact.
Now assume that the space $G_{\vec S}$ is countably pracompact and Baire. By Proposition~\ref{G_S*CPCompact}, the subgroup $S-S$ has finite index in $G$ and  $S\subseteq C-S$ for some countable set $C\subseteq S$. By Proposition~\ref{G_S*Baire}, $S=S-S$ and by Proposition~\ref{G_S*Compact},  the space $G_{\vec S}$ is compact.
\end{proof}
A space $X$ is called {\it a $P$-space} if every $G_\delta$-subset of $X$ is open.
\begin{proposition}\label{G_S*PSp} The space $G_{\vec S}$ is a $P$-space if and only if each countable subset $C\subseteq S$ is contained in the set $a-S$ for some element $a\in S$.
\end{proposition}
\begin{proof} The characterization is trivial if $S=\{0\}$. In this case the space $S_{\vec S}$ is discrete and hence a $P$-space. Also $S\subseteq -S$. So, we assume that $S\ne\{0\}$.
Given any countable subset $C\subseteq S$, consider the $G_\delta$-set $W=\bigcap_{c\in C}(\{0\}\cup(c+S))$ in $G_{\vec S}$. If $G_{\vec S}$ is a $P$-space, then $W$ is a neighborhood of zero. Then there exists an element $b\in S$ such that $b+S\subseteq W$. If $b\ne 0$, then put $a:=b$. If $b=0$, then take any point $a\in S\setminus\{0\}$ and conclude that $a+S\subseteq S=b+S\subseteq W$.  In both cases there exists an element $a\in S\setminus\{0\}$ such that $a+S\subseteq W$. Then for every $c\in C$ the inclusion $0\ne  a\in a+S\subseteq\{0\}\cup(c+S)$ implies $a\in c+S$ and hence $c\in a-S$ and finally $C\subseteq a-S$.
\smallskip
Now assuming that each countable subset of $S$ is contained in the set $a-S$ for some $a\in S$, we shall prove that $G_{\vec S}$ is a $P$-space. Given any point $x\in G$ and a sequence $(U_n)_{n\in\omega}$ of neighborhoods of $x$, we should prove that the intersection $W:=\bigcap_{n\in\omega}U_n$ is a neighborhood of $x$. For every $n\in\omega$ choose a point $c_n\in S$ such that
$x+c_n+S\subseteq U_n$. By our assumption, there exists an element $a\in S$ such that $\{c_n:n\in\omega\}\subset
a-S$. Then $x+a+S\subseteq x+c_n+S\subseteq U_n$ for each $n\in\omega$. Hence
$x+a+S\subseteq \bigcap_{n\in\omega}U_n=W$ and thus $W$ is a neighborhood of $x$ in $G_{\vec S}$.
\end{proof}
Propositions~\ref{G_S*2PCompact}, \ref{G_S*PCompact}, \ref{G_S*PSp} imply the following topological characterization of the $2$-pseudocompact\-ness of the paratopological group $G_{\vec S}$.
\begin{corollary}\label{G_S*2PCompactTop} The paratopological group $G_{\vec S}$ is $2$-pseudocompact if  and only if $G_{\vec S}$ is a feebly compact P-space.
\end{corollary}
On the other hand, Corollaries~\ref{G_S*CompactTop} and \ref{G_S*2PCompactTop} imply another characterization of the compactness of  $G_{\vec S}$.
\begin{corollary} The space $G_{\vec S}$ is compact if and only if $G_{\vec S}$
is a countably pracompact P-space.
\end{corollary}
The following example shows that a counterpart of Proposition~\ref{G_SFCompact} does not hold for the paratopological groups $G_{\vec S}$.
\begin{example}\label{G_S*FCompact} For the group $G=\IR$ and the subsemigroup $S:=[0,\infty)$ of $G$,
the paratopological group $G_{\vec S}$ is $\omega$-precompact and not Lindel\"of. By
Propositions~\ref{G_S*CPCompact} and \ref{G_S*2PCompact}, the paratopological group $G_{\vec S}$ is
countably pracompact but not $2$-pseudocompact.
\end{example}
\begin{example}\label{Zomega1} There are an Abelian group $G$ and a subsemigroup $S$ of the group $G$
such that the paratopological groups $G_S$ and $G_{\vec S}$ have the following properties:
\begin{enumerate}
\item $G_S$ is $T_0$, sequentially compact, not totally countably compact,
not precompact, and not a topological group;
\item $G_{\vec S}$ is $T_1$, $2$-pseudocompact, not
countably pracompact, not precompact, and not a topological group.
\end{enumerate}
\end{example}
\begin{proof} Let $G=\bigoplus\limits_{\alpha\in\omega_1} \IZ$ be the direct sum of $\omega_1$ many copies of the group $\IZ$. In the group $G$ consider the subsemigroup $$S=\{0\}\cup\big\{(x_\alpha)_{\alpha\in\omega_1}\in G:\exists \beta\in\omega_1\; (x_\beta>0\mbox{ and } \forall \alpha>\beta\;\;x_\alpha=0)\big\}.$$
Since $S\cap (-S)=\{0\}$, we can apply Propositions~\ref{G_ST0}, \ref{G_S*T} and conclude that the space $G_S$ is $T_0$ and
 the space $G_{\vec S}$ is $T_1$.
Since $G=S-S$ and for each countable infinite subset $C$ of $S$ there is an element $a\in S$ such that $C\subset
a-S$, we see that
by Proposition~\ref{G_S2PCompact}, the group $G_S$ is sequentially compact and
by Proposition~\ref{G_S*2PCompact}, the group $G_{\vec S}$ is $2$-pseudocompact.
By Proposition~\ref{G_S*CPCompact}, the group $G_{\vec S}$ is not countably pracompact.
By Proposition~\ref{G_SCompact}, the group $G_S$ is not totally countably compact and not
precompact. Therefore the group $G_{\vec S}$ is not precompact, too.
\end{proof}
\begin{remark} In \cite{San}, S\'anchez proved that the group $G_{\vec S}$ from Example~\ref{Zomega1} is a $P$-space.
\end{remark}
The following proposition was announced in \cite{RavRez}.
\begin{proposition}\label{T1CCSquareTop}
Each $T_1$ paratopological group $G$ with a countably compact square $G\times G$ is a topological
group.
\end{proposition}
\begin{proof}
Indeed, $G^*=\{(x,x^{-1})\in G\times G:x\in G\}$ is a closed (and hence a
countably compact) subgroup of $G\times G$. Since $G^*$ has a base at the identity consisting
of symmetric open neighborhoods, it is a topological group. By Lemma \ref{PseudoTGPre},
the group $G^*$ is precompact. Since a projection $\pi_1: G\times G\to F$ onto the first coordinate is
a continuous surjective homomorphism and $\pi_1(G^*)=G$, the group $G$ is precompact too.
By Proposition~\ref{PseTG}, $G$ is a topological group.
\end{proof}
Now we see that Proposition~\ref{T1CCSquareTop}
cannot be generalized to $T_1$ sequentially pracompact groups
by Example~\ref{BanakhEx} and
to $T_0$ sequentially compact groups by Example~\ref{Zomega1}.
\subsection{Cowide topologies}\label{Cowide}
Topologies $\tau$ and $\sigma$ on a set $X$ are defined to be {\it cowide} if for any
nonempty sets $U\in\tau$ and $V\in\sigma$ the intersection $U\cap V$ is nonempty. If a
topology $\sigma$ is cowide to itself then  the topology $\sigma$ will be called {\it wide}.
For two topologies  $\tau$ and $\sigma$ on a set $X$ by $\tau\vee\sigma$ we denote their supremum; it is generated by the base $\{U\cap V:U\in\tau,
V\in\sigma\}$.
\begin{lemma}\label{ClosureOfOpen} Let $\tau$ and $\sigma$ be cowide topologies on a set $X$.
Then $\ol{W}^{\tau\vee\sigma}=\ol{W}^\tau$ for each set $W\in\tau$. Moreover,
if the topology $\sigma$ is wide, then
$\ol{W}^{\tau\vee\sigma}=\ol{W}^\tau$ for each set $W\in\tau\vee\sigma$.
\end{lemma}
\begin{proof} The inclusion $\tau\vee\sigma\supseteq\tau$ implies the inclusion
$\ol{W}^{\tau\vee\sigma}\subseteq \ol{W}^\tau$ for any subset $W\subseteq X$.
Suppose that $W\in\tau$, $x\in\ol{W}^\tau$
and $U\cap V$ is an arbitrary neighborhood of the point $x$ such that $U\in\tau$ and $V\in\sigma$.
Then $W\cap U$ is a nonempty $\tau$-open set. Since the topologies
$\tau$ and $\sigma$ are cowide, the intersection $W\cap U\cap V$ is nonempty. Therefore $x\in \ol{W}^{\tau\vee\sigma}$ and hence $\ol{W}^\tau=\ol{W}^{\tau\vee\sigma}$.
Now suppose that the topology $\sigma$ is wide, $W\in\tau\vee\sigma$, $x\in\ol{W}^\tau$
and $U\cap V\in\tau\vee\sigma$ is any neighborhood of the point $x$ such that $U\in\tau$ and $V\in\sigma$.
The intersection $W\cap U$ is nonempty and belongs to the topology $\tau\vee\sigma$. Therefore there exist sets $U'\in\tau$ and $V'\in\sigma$ such that
$\0\ne U'\cap V'\subseteq W\cap U$.
Since the topology $\sigma$ is wide, the intersection $V\cap V'$ is nonempty.
Since the topologies
$\tau$ and $\sigma$ are cowide, the intersection $U'\cap V\cap V'\subseteq W\cap U\cap V$
is nonempty. Therefore $x\in \ol{W}^{\tau\vee\sigma}$.
\end{proof}
\subsubsection{Semiregularizations}
\begin{lemma}\label{CowideRegularization} Let $\tau$ and $\sigma$ be cowide topologies on a set $X$ such that the topology $\sigma$ is wide. Then $(\tau\vee\sigma)_{sr}=\tau_{sr}$.
\end{lemma}
\begin{proof} The inclusion $\tau\vee\sigma\supseteq\tau$ implies the inclusion
$(\tau\vee\sigma)_{sr}\supseteq\tau_{sr}$. So we show the opposite inclusion.
Let $W=\int_{\tau\vee\sigma}\ol{W}^{\tau\vee\sigma}$
be any regular open set in the topology $\tau\vee\sigma$. Lemma~\ref{ClosureOfOpen} implies that
$\ol{W}^{\tau\vee\sigma}=\ol{W}^\tau$.
The inclusion $\tau\vee\sigma\supseteq\tau$ implies the inclusion
$\int_{\tau\vee\sigma}\ol{W}^\tau\supseteq\int_\tau\ol{W}^\tau$.
So we show the opposite inclusion. Let $x\in \int_{\tau\vee\sigma}\ol{W}^\tau$
be any point.
Therefore there exist sets $U'\in\tau$ and $V'\in\sigma$ such that
$x\in U'\cap V'\subseteq \ol{W}^\tau$. Then $\ol{U'\cap V'}^\tau\subseteq \ol{W}^\tau$.
Since the topologies $\tau$ and $\sigma$ are cowide, the set $U'\cap V'$ is $\tau$-dense in the set $U'$.
Therefore $\ol{U'\cap V'}^\tau=\ol{U'}^\tau$. Hence
$x\in U'\subseteq \ol{U'}^\tau\subseteq \ol{W}^\tau$ and thus $x\in \int_\tau\ol{W}^\tau$.
\end{proof}
Remark that similar constructions were considered by Bourbaki in order to
provide non-regular topologies finer than regular, see Exercise 20 for Chapter 1, $\S$8
in~\cite{Bou2}.
\subsubsection{Feeble compactness}
Lemmas~\ref{CowideRegularization} and~\ref{RegPCompact} imply the following
\begin{proposition}\label{SupPseudo} Let $\tau$ and $\sigma$ be cowide topologies on a set $X$ such that
the space $(X,\tau)$ is feebly compact and the topology $\sigma$ is wide. Then the
space $(X,\tau\vee\sigma)$ is feebly compact, too.
\end{proposition}

{\subsubsection{$H$-compactness} We recall that a space $X$ is {\it $H$-compact} iff for each open cover $\mathcal U$ of the space $X$ there exists a finite subfamily $\mathcal V\subseteq\mathcal U$ such that $X=\bigcup_{V\in\mathcal V} \ol V$.
It is easy to see that a regular $H$-compact space is compact and
a space is $H$-compact if and only if its semiregularization is $H$-compact. This characterization combined with Lemma~\ref{CowideRegularization} imply the following proposition.
\begin{proposition}\label{SupHclosed} Let $\tau$ and $\sigma$ be cowide topologies on a set $X$ such that
the space $(X,\tau)$ is $H$-compact and the topology $\sigma$ is wide. Then the
space $(X,\tau\vee\sigma)$ is $H$-compact, too.
\end{proposition}
\begin{proposition}\label{p:3.30} For any paratopological group $(G,\tau)$ the following conditions are equivalent:
\begin{enumerate}
\item The space $(G,\tau)$ is $H$-compact.
\item The space $(G,\tau_{sr})$ is $H$-compact.
\item The space $(G,\tau_{sr})$ is compact.
\end{enumerate}
If the space $(G,\tau)$ is quasiregular, then the conditions \textup{(1)--(3)} are equivalent to
\begin{itemize}
\item[\textup{(4)}] The space $(G,\tau)$ is compact.
\end{itemize}
\end{proposition}
\begin{proof} The implications $(1)\Leftrightarrow(2)\Leftarrow(3)\Leftarrow(4)$ are trivial and
hold for all spaces. To prove that $(2)\Ra(3)$, assume that the space $(G,\tau_{sr})$ is $H$-compact.
By \cite[Ex. 1.9]{Rav2}, \cite[p. 31]{Rav3}, and \cite[p. 28]{Rav3}, $(G,\tau_{sr})$ is a regular
paratopological group, so it is compact.
Now assuming that the space $(G,\tau)$ is quasiregular, we shall prove that $(3)\Ra(4)$. Assume that the space $(G,\tau_{sr})$ is compact. By Proposition~\ref{PseTG}, the quasiregular compact paratopological group $(G,\tau_{sr})$ is a topological group. Then by Lemma~\ref{RegTG},
$(G,\tau)$ is a topological group, too. Therefore the space $(G,\tau)$ is regular and
hence  $\tau_{sr}=\tau$ and the space $(G,\tau)=(G,\tau_{sr})$ is compact.
\end{proof}
}
\subsubsection{Applications of cowide topologies to the cone topologies} From now on and till the end of the section we assume that $G$ is an Abelian group and $S\ni 0$ is a subsemigroup in $G$. Let $\tau$ be any semigroup topology on $G$. It is easy to see that the cone topology $\tau_S$ defined at the beginning of Section~\ref{sec:cone-top} coincides with the supremum of the topologies of the spaces $(G,\tau)$ and $G_S$.
Similarly, the shift cone topology $\tau_{\vec S}$ is the superemum of topologies of the spaces $(G,\tau)$ and $G_{\vec S}$.
It is clear that the group $G$ endowed with the topology $\tau_{\vec S}$ is a paratopological group.
\begin{proposition}\label{WideS} The following conditions are equivalent:
\begin{enumerate}
\item The topology of the space $G_S$ is wide.
\item The topology of the space $G_{\vec S}$ is wide.
\item $G=S-S$.
\end{enumerate}
\end{proposition}
\begin{proof}
Since the topology of the space $G_S$ is coarser than
the topology of the space $G_{\vec S}$, we have implication $(2)\Rightarrow(1)$.
\smallskip
$(1)\Rightarrow(3)$ Let $x\in G$ be any element. If the topology of the space $G_S$ is wide, then
the intersection $S\cap (x+S)$ of open subsets $S$ and $x+S$ of $G_S$ is nonempty and hence $x\in
S-S$.
\smallskip
$(3)\Rightarrow(2)$ Observe that every nonempty open subset of the group $G_{\vec S}$ contains the
set $x+S$ for some element $x\in G$. Let $U_1$, $U_2$ be any nonempty open subsets of the space
$G_{\vec S}$. Then there are elements $x_1,x_2\in G$ such that $x_1+S\subseteq U_1$ and $x_2+S\subset
U_2$. If $G=S-S$, then there are elements $y_1,y_2,z_1,z_2\in S$ such that $x_1=y_1-z_1$ and
$x_2=y_2-z_2$. Then the intersection $U_1\cap U_2\supset(x_1+S)\cap (x_2+S)$ contains the element
$y_1+y_2$ and hence is not empty. This means that the topology of the space $G_{\vec S}$ is wide.
\end{proof}

\begin{proposition}\label{CowideS} The following conditions are equivalent:
\begin{enumerate}
\item The topology of the space $G_S$ and the topology $\tau$ are cowide.
\item The topology of the space $G_{\vec S}$ and the topology $\tau$ are cowide.
\item $S$ is a dense subsemigroup of the group $(G,\tau)$.
\end{enumerate}
\end{proposition}
\begin{proof} Since the topology of the space $G_S$ is coarser than
the topology of the space $G_{\vec S}$, we have implication $(2)\Rightarrow(1)$.
The implication $(1)\Rightarrow(3)$ follows from the opennes of the set $S$ is the paratopological group $G_{\vec S}$.
The implication $(3)\Rightarrow(2)$ follows from the observation that
each nonempty open subset of the space $G_{\vec S}$ contains a set $x+S$, $x\in S$, which is dense the space $(G,\tau)$.
\end{proof}
\begin{proposition} Assume that the semigroup $S$ is dense in $(G,\tau)$ and $G=S-S$.
Then the following conditions are equivalent:
\begin{enumerate}
\item The space $(G,\tau)$ is feebly compact.
\item  The space $(G,\tau_S)$ is feebly compact.
\item The space $(G,\tau_{\vec S})$ is feebly compact.
\end{enumerate}
\end{proposition}
\begin{proof} The implications $(3)\Rightarrow(2)\Rightarrow(1)$ are obvious, and  $(1)\Rightarrow(3)$
 follows from Propositions~\ref{CowideS}, \ref{WideS}, and \ref{SupHclosed}.
\end{proof}
\begin{proposition} Assume that the semigroup $S$ is dense in $(G,\tau)$ and $G=S-S$.
Then the following conditions are equivalent:
\begin{enumerate}
\item The space $(G,\tau)$ is $H$-compact.
\item The space $(G,\tau_S)$ is $H$-compact.
\item The space $(G,\tau_{\vec S})$ is $H$-compact.
\end{enumerate}
\end{proposition}
\begin{proof}
The implications $(3)\Rightarrow(2)\RT(1)$ are obvious, and $(1)\Rightarrow(3)$ follows from Propositions~\ref{CowideS}, \ref{WideS}, and \ref{SupHclosed}.
\end{proof}

\section{Acknowledgements}
The second author is grateful to Artur Tomita and Manuel Sanchis for providing
him their articles.
His sincere thanks also go to Minnan Normal University (Fujian,
China) and to University of Hradec Kr\'alov\'e (Czech Republic)
for the hospitality and opportunity
 to give lectures on the topic of the
present paper, to Mikhail Tkachenko who has
inspired him to prove Proposition~\ref{2PseP}, and
to {\tt MathOverflow} user Dog$\_$69, who pointed him\footnote{\tt https://math.stackexchange.com/q/3138278} a reference
from~\cite{Bou2} related to Lemma~\ref{CowideRegularization}.

\end{document}